\documentclass[10pt,a4paper]{amsart}

\usepackage[austrian,english]{babel}
\usepackage[latin1]{inputenc}
\usepackage[a4paper]{geometry}
\usepackage{color}
\usepackage{tikz}
\usetikzlibrary{arrows,decorations.pathreplacing,patterns,shapes}
\usepackage[font=footnotesize]{caption}
\usepackage{enumerate}
\usepackage{xparse}

\usepackage{amsmath}
\usepackage{amssymb}
\usepackage{faktor}
\usepackage{mathtools}

\usepackage{hyperref}

\newcommand{\N}{\mathbb{N}}

\newcommand{\coQ}{\check{Q}}
\newcommand{\R}{\mathbb{R}}
\newcommand{\W}{\widetilde{W}}
\newcommand{\Z}{\mathbb{Z}}
\renewcommand{\S}{\mathfrak{S}}

\newcommand{\act}{\cdot}
\newcommand{\abs}[1]{\left| #1 \right|}
\newcommand{\ac}{A_{\circ}}
\newcommand{\floor}[1]{\left\lfloor #1 \right\rfloor}

\newcommand{\skal}[1]{\langle #1 \rangle}
\newcommand{\op}[1]{\operatorname{#1}}

\newcommand{\affS}{{\smash{\stackrel{\sim}{\smash{\mathfrak{S}}\rule{0pt}{1.05ex}}}}}

\newcommand{\sk}{

\smallskip}

\NewDocumentEnvironment{mymathenvironment}{m m}{\begin{#1}\label{#1:#2}}{\end{#1}}

\NewDocumentEnvironment{mythrm}{m}{\begin{mymathenvironment}{Theorem}{#1}}{\end{mymathenvironment}}
\NewDocumentEnvironment{mycor}{m}{\begin{mymathenvironment}{Corollary}{#1}}{\end{mymathenvironment}}
\NewDocumentEnvironment{myconj}{m}{\begin{mymathenvironment}{Conjecture}{#1}}{\end{mymathenvironment}}
\NewDocumentEnvironment{mylem}{m}{\begin{mymathenvironment}{Lemma}{#1}}{\end{mymathenvironment}}
\NewDocumentEnvironment{myprop}{m}{\begin{mymathenvironment}{Proposition}{#1}}{\end{mymathenvironment}}

\newcommand{\refx}[2]{#2~\ref{#2:#1}}

\newcommand{\refcj}[1]{\refx{#1}{Conjecture}}
\newcommand{\refc}[1]{\refx{#1}{Corollary}}
\newcommand{\reff}[1]{\refx{#1}{Figure}}
\newcommand{\refl}[1]{\refx{#1}{Lemma}}
\newcommand{\refp}[1]{\refx{#1}{Proposition}}
\newcommand{\refq}[1]{\eqref{eq:#1}}
\newcommand{\refs}[1]{\refx{#1}{Section}}
\newcommand{\reft}[1]{\refx{#1}{Theorem}}

\newcommand{\refi}[1]{\eqref{item:#1}}

\newcommand{\drawrunner}[2]{
\pgfmathsetmacro{\runner}{#1}
\pgfmathsetmacro{\bead}{#2}
\pgfmathsetmacro{\b}{-\bead}
\foreach \level in {\b,...,\min}{
\pgfmathtruncatemacro{\z}{\runner - \n*\level}
\draw(\runner,\level) node[shape=circle,draw,inner sep=.1pt,minimum size=1.5em]{\footnotesize{\z}};
\draw[yshift=1cm](\runner,\min)node{$\vdots$};
}
\pgfmathsetmacro{\b}{-1-\bead}
\foreach \level in {\max,...,\b}{
\pgfmathtruncatemacro{\z}{\runner - \n*\level}
\draw(\runner,\level) node{\footnotesize{\z}};
}
\draw[yshift=-.7cm](\runner,\max)node{$\vdots$};
}

\begin{document}

\title{Rational Shi tableaux and the skew length statistic}
\author{Robin Sulzgruber}\thanks{Research supported by the Austrian Science Fund (FWF), grant S50-N15 in the framework of the Special Research Program ``Algorithmic and Enumerative Combinatorics'' (SFB F50).}
\date{September~2016}
\address{Fakult{\"a}t f{\"u}r Mathematik, Universit{\"a}t Wien, Oskar-Morgenstern-Platz 1, 1090 Wien, Austria}
\email{robin.sulzgruber@univie.ac.at}
\begin{abstract}
We define two refinements of the skew length statistic on simultaneous core partitions.
The first one relies on hook lengths and is used to prove a refined version of the theorem stating that the skew length is invariant under conjugation of the core.
The second one is equivalent to a generalisation of Shi tableaux to the rational level of Catalan combinatorics.
These rational Shi tableaux encode dominant $p$-stable elements in the affine symmetric group.
We prove that the rational Shi tableau is injective, that is, each dominant $p$-stable affine permutation is determined uniquely by its Shi tableau.
Moreover, we provide a uniform generalisation of rational Shi tableaux to Weyl groups, and conjecture injectivity in the general case.

\end{abstract}

\thispagestyle{empty}
\maketitle




\pagenumbering{arabic}
\pagestyle{headings}

The skew length is a statistic on simultaneous core partitions, which counts the number of certain cells in the Young diagram of the core. It was invented by Armstrong~\cite{AHJ2014} who used it to give a combinatorial formula for the rational $q,t$-Catalan numbers (see also~\cite{GM2013,ALW2014:sweep_maps})
\begin{align*}
C_{n,p}(q,t)
=\sum_{\kappa\in\mathfrak{C}_{n,p}} q^{\ell(\kappa)} t^{(n-1)(p-1)/2-\op{skl}(\kappa)}.
\end{align*}
Here $n$ and $p$ are relatively prime, $\mathfrak C_{n,p}$ denotes the set of $n,p$-cores, $\ell(\kappa)$ is the length of $\kappa$ and $\op{skl}(\kappa)$ denotes the skew length (see \refs{conj} for all definitions). 
The intriguing property of these polynomials is the apparent symmetry $C_{n,p}(q,t)=C_{n,p}(t,q)$, for which there is no proof except in very special cases.

The skew length was further studied in~\cite{Xin2015} and ~\cite{CDH2015}. For coprime $n$ and $p$ it was proven that $\op{skl}(\kappa)$ is invariant under conjugation of the core and independent of whether $\kappa$ is viewed as an $n,p$-core or as a $p,n$-core.
The first main contribution of our paper is a refinement of the skew length statistic in terms of the multiset of hook lengths of certain cells of the core.
We obtain a refined version of the two results mentioned above in \reft{Hab} providing a conceptual explanation for these phenomena.
Moreover our method of proof relies on induction and does not use rational Dyck paths.
This has the advantage that we do not need to assume that $n$ and $p$ are relatively prime.
\reft{Hab} therefore extends to previously untreated cases.

\smallskip
Shi tableaux~\cite{FTV2011} encode the dominant regions of the ($m$-extended) Shi arrangement, which are counted by the Fu{\ss}--Catalan numbers. The Shi arrangement and its regions are intimately related to the affine symmetric group. In~\cite{GMV2014} Gorsky, Mazin and Vazirani define the so called $p$-stable affine permutations, which can be seen as a generalisation of the regions of the Shi arrangement to the rational level of Catalan combinatorics. In \refs{shitab} we define rational Shi tableaux, which encode dominant $p$-stable affine permutations. Our main result concerning rational Shi tableaux is \reft{shiA}, asserting that each dominant $p$-stable affine permutation is uniquely determined by its Shi tableau.
\reft{shiA} is especially interesting in view of the fact that the rational Pak--Stanley labelling~\cite{GMV2014} can be obtained by taking the row-sums of the rational Shi tableau.

\smallskip
The Shi arrangement can be defined for any irreducible crystallographic root system $\Phi$.
Furthermore, $p$-stable affine Weyl group elements have been considered in~\cite{Thiel2015}.
We provide a uniform definition of the rational Shi tableaux of dominant $p$-stable elements of the affine Weyl group of $\Phi$.
We conjecture that dominant $p$-stable Weyl group elements are determined uniquely by their Shi tableaux
(\refcj{shi}).

\smallskip
In \refs{codinv} we tie our previous results together by relating dominant $p$-stable affine permutations to rational Dyck paths and simultaneous cores. This is achieved using the Anderson bijection of~\cite{GMV2014}. We show that the rational Shi tableau of a dominant $p$-stable affine permutation can be computed from the corresponding rational Dyck path as the codinv tableau of the path (\reft{dinvshi}). The codinv statistic on rational Dyck paths, that is, the sum of the entries of the codinv tableau, corresponds naturally to the skew length statistic on cores (\refc{codinv}). Hence we find that rational Shi tableaux can in fact be regarded as another refinement of the skew length statistic. The contents of \refs{codinv} and in particular the connection between rational Shi tableaux and the zeta map on rational Dyck paths~\cite{GM2013,ALW2014:sweep_maps} (\reft{zeta}) are the main tools used in the proof of \reft{shiA}.

\smallskip
An extended abstract~\cite{Sul2016} of this paper has appeared in the proceedings of FPSAC~2016 in Vancouver.

\section{Skew-length and conjugation}\label{Section:conj}

\begin{figure}[t]
\begin{center}
\begin{tikzpicture}[scale=.7]
\begin{scope}
\draw[very thick, green, fill=green!50, rounded corners=3mm]
(0,13) rectangle (1,10)
(0,6) rectangle (1,9)
(0,1) rectangle (1,5)
(1,13) rectangle (2,10)
(1,6) rectangle (2,9)
(1,1) rectangle (2,5)
(4,13) rectangle (5,10)
(4,7) rectangle (5,9)
(8,11) rectangle (9,13)
(9,11) rectangle (10,13)
(11,11) rectangle (12,13)
;
\draw[very thick, red, fill=red, opacity=.4, rounded corners=3mm]
(0,12) rectangle (6,13)
(8,12) rectangle (12,13)
(14,12) rectangle (15,13)
(16,12) rectangle (18,13)
(0,9) rectangle (6,10)
(0,8) rectangle (5,9)
(0,7) rectangle (5,8)
(0,2) rectangle (2,3)
;
\draw[very thick] (0,0)--(1,0)--(1,1)--(2,1)--(2,5)--(3,5)--(3,6)--(4,6)--(4,7)--(5,7)--(5,9)--(6,9)--(6,10)--(7,10)--(7,11)--(13,11)--(13,12)--(19,12)--(19,13)--(0,13)--cycle;
\draw (1,1)--(1,13)
(2,5)--(2,13)
(3,6)--(3,13)
(4,7)--(4,13)
(5,9)--(5,13)
(6,10)--(6,13)
(7,11)--(7,13)
(8,11)--(8,13)
(9,11)--(9,13)
(10,11)--(10,13)
(11,11)--(11,13)
(12,11)--(12,13)
(13,12)--(13,13)
(14,12)--(14,13)
(15,12)--(15,13)
(16,12)--(16,13)
(17,12)--(17,13)
(18,12)--(18,13)
(0,1)--(1,1)
(0,2)--(2,2)
(0,3)--(2,3)
(0,4)--(2,4)
(0,5)--(2,5)
(0,6)--(3,6)
(0,7)--(4,7)
(0,8)--(5,8)
(0,9)--(5,9)
(0,10)--(6,10)
(0,11)--(7,11)
(0,12)--(13,12);
\draw[xshift=5mm,yshift=5mm]
(0,12) node{31}
(0,11) node{24}
(0,10) node{17}
(0,9) node{15}
(0,8) node{13}
(0,7) node{12}
(0,6) node{10}
(0,5) node{8}
(0,4) node{6}
(0,3) node{5}
(0,2) node{4}
(0,1) node{3}
(0,0) node{1}
(1,12) node{29}
(1,11) node{22}
(1,10) node{15}
(1,9) node{13}
(1,8) node{11}
(1,7) node{10}
(1,6) node{8}
(1,5) node{6}
(1,4) node{4}
(1,3) node{3}
(1,2) node{2}
(1,1) node{1}
(2,12) node{24}
(2,11) node{17}
(2,10) node{10}
(2,9) node{8}
(2,8) node{6}
(2,7) node{5}
(2,6) node{3}
(2,5) node{1}
(3,12) node{22}
(3,11) node{15}
(3,10) node{8}
(3,9) node{6}
(3,8) node{4}
(3,7) node{3}
(3,6) node{1}
(4,12) node{20}
(4,11) node{13}
(4,10) node{6}
(4,9) node{4}
(4,8) node{2}
(4,7) node{1}
(5,12) node{17}
(5,11) node{10}
(5,10) node{3}
(5,9) node{1}
(6,12) node{15}
(6,11) node{8}
(6,10) node{1}
(12,11) node{1}
(11,11) node{2}
(10,11) node{3}
(9,11) node{4}
(8,11) node{5}
(7,11) node{6}
(18,12) node{1}
(17,12) node{2}
(16,12) node{3}
(15,12) node{4}
(14,12) node{5}
(13,12) node{6}
(12,12) node{8}
(11,12) node{9}
(10,12) node{10}
(9,12) node{11}
(8,12) node{12}
(7,12) node{13}
;
\end{scope}
\end{tikzpicture}
\caption{A simultaneous core $\kappa\in\mathfrak C_{7,16}$ with the multisets $H_{7,16}(\kappa)$ in red and $H_{16,7}(\kappa)$ in green.}
\label{Figure:core}
\end{center}
\end{figure}
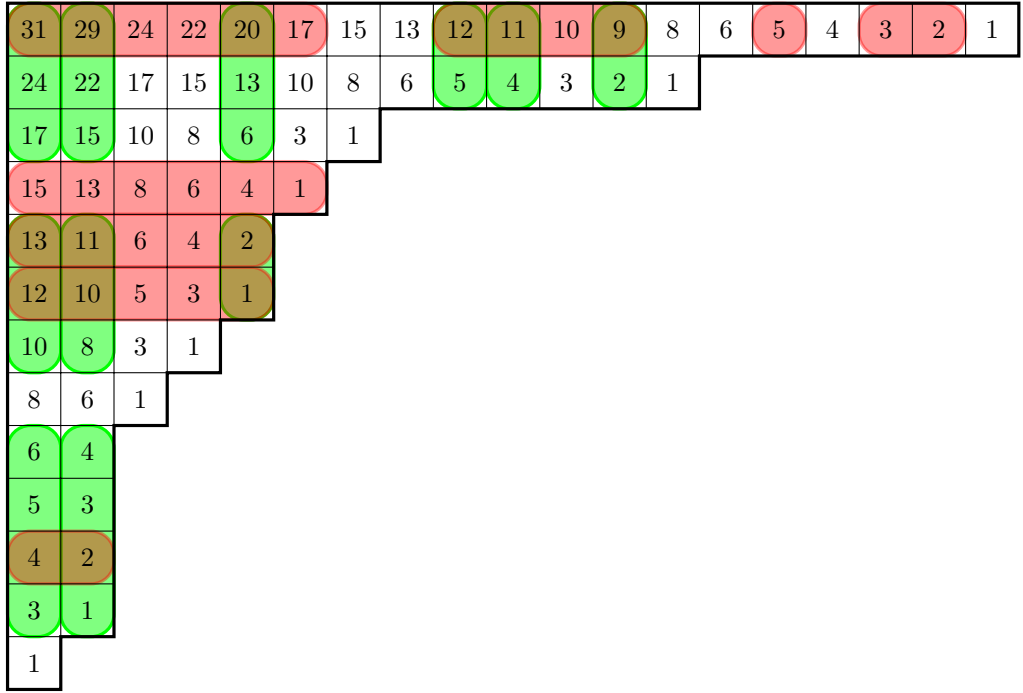

Throughout this section $n$ and $p$ are positive integers (not necessarily coprime!). For any positive integer $z$ we set $[z]=\{1,\dots,z\}$.
A \emph{partition} $\lambda$ is a weakly decreasing sequence $\lambda_1\geq\lambda_2\geq\dots\geq\lambda_{\ell}>0$ of positive integers.
The number of \emph{parts} or \emph{summands} $\lambda_i$ is called the \emph{length} of the partition and is denoted by $\ell(\lambda)$. 
The sum $\sum_i{\lambda_i}$ is called the \emph{size} of $\lambda$. 
A partition is often identified with its \emph{Young diagram} $\{(i,j):i\in[\ell(\lambda)],j\in[\lambda_i]\}$.
We call the elements of the Young diagram \emph{cells} of the partition.
The \emph{conjugate} of the partition $\lambda$ is the partition $\lambda'$ where $\lambda'_i=\#\{j:\lambda_j\geq i\}$.
Equivalently, if we view partitions as Young diagrams then $(j,i)\in\lambda'$ if and only if $(i,j)\in\lambda$.
The \emph{hook length} of a cell of $\lambda$ is defined as $h_{\lambda}(i,j)=\lambda_i-j+\lambda'_j-i+1$.
A partition $\kappa$ is called \emph{$n$-core} if no cell of $\kappa$ has hook length $n$.
A partition is called $n,p$-core if it is both an $n$-core and a $p$-core.
We denote the set of $n$-core partitions by $\mathfrak C_n$ and the set of simultaneous $n,p$-cores by $\mathfrak C_{n,p}$.
Clearly $\lambda\in\mathfrak{C}_n$ if and only if $\lambda'\in\mathfrak{C}_n$.

\smallskip
Let $\lambda$ be a partition.
Given the hook length $h$ of a cell $x$ in the top row of $\lambda$ we denote by $H^c(h)$ the set of hook lengths of cells in the same column as $x$. 
Given the hook length $h$ of a cell $x$ in the first column of $\lambda$ we denote by $H^r(h)$ the set of hook lengths of cells in the same row as $x$. 
If $h$ is not the hook length of a suitable cell, set $H^r(h)=\emptyset$ resp.~$H^c(h)=\emptyset$.
For example, in \reff{core} we have $H^c(12)=\{12,5\}$ and $H^r(12)=\{12,10,5,3,1\}$ and $H^r(11)=\emptyset$.

\smallskip
A subset $A\subseteq\Z$ is called \emph{abacus} if there exist integers $a,b\in\Z$ such that $z\in A$ for all $z$ with $z<a$, and $z\notin A$ for all $z$ with $z>b$.
We call the elements of $A$ \emph{beads} and the elements of $\Z-A$ \emph{gaps}.
An abacus is \emph{normalised} if zero is a gap and there are no negative gaps.
An abacus is \emph{balanced} if the number of positive beads equals the number of non-positive gaps.
An abacus $A$ is \emph{$n$-flush} if $z-n\in A$ for all $z\in A$.
The theorem below is a version of the classical result that $n$-cores correspond to abacus diagrams that are $n$-flush (see for example \cite[2.7.13]{JamesKerber}).
Let $\kappa\in\mathfrak C_n$ be an $n$-core with maximal hook length $m$, that is, $m$ is the hook length of the top left corner.
Define
$\alpha(\kappa)=\{\kappa_i+1-i:i\geq1\}$ and 
$\beta(\kappa)=H^c(m)\cup\{z\in\Z:z<0\}$,
where by convention $\kappa_i=0$ when $i>\ell(\kappa)$. Note that $\beta(\kappa)=\{z+\ell(\kappa)-1:z\in\alpha(\kappa)\}$.

\begin{mythrm}{coreabacus} The map $\alpha$ is a bijection between $n$-cores and balanced $n$-flush abaci. The map $\beta$ is a bijection between $n$-cores and normalised $n$-flush abaci.
\end{mythrm}

For our purposes in this section we need the following simple consequence of \reft{coreabacus}.
 
\begin{mylem}{flush} Let $\kappa\in\mathfrak C_n$ be an $n$-core
and $z\geq0$. Then $z+n\in H^r(h)$ implies $z\in H^r(h)$, and $z+n\in H^c(h)$ implies $z\in H^c(h)$.
\end{mylem}

In fact, the set of $n$-cores is characterised by the property that $z+n\in H^c(m)$ implies $z\in H^c(m)$ for all $z\geq0$, where $m$ is the maximal hook length of $\kappa$.

\smallskip
Let $\kappa\in\mathfrak C_{n,p}$ be a simultaneous core with maximal hook length $m$. Moreover choose an element $h\in H^c(m)$. 
We call the row of $\kappa$ whose leftmost cell has hook length $h$ an $n$-row (resp.~$p$-row) if $h+n\notin H^c(m)$ (resp.~$h+p\notin H^c(m)$). 
Similarly, given $h\in H^r(m)$ we call the column whose top cell has hook length $h$ an $n$-column (resp.~$p$-column) if $h+n\notin H^r(m)$ (resp.~$h+p\notin H^r(m)$). 
Denote by $H_{n,p}(\kappa)$ the multiset of hook lengths of cells that are contained both in an $n$-row and in a $p$-column.
In \reff{core} the leftmost hook lengths of the $7$-rows are $31,15,13,12$ and $4$, and the top hook lengths of the $7$-columns are $31,29,20,12,11$ and $9$.

\smallskip
This section's main result is a surprising symmetry property of the multiset $H_{n,p}(\kappa)$.

\begin{mythrm}{Hab} Let $\kappa\in\mathfrak C_{n,p}$ be an $n,p$-core. Then $H_{n,p}(\kappa)=H_{p,n}(\kappa)$.
\end{mythrm}

\begin{proof} We prove the claim by induction on the size of $\kappa$. 
Denote by $\bar\kappa$ the partition obtained from $\kappa$ by deleting the first column.
Clearly $\bar\kappa\in\mathfrak C_{n,p}$ and we may assume that $H_{n,p}(\bar\kappa)=H_{p,n}(\bar\kappa)$.

Let $m$ denote the maximal hook length in $\kappa$.
Note that each $n$-row of $\bar\kappa$ is an $n$-row of $\kappa$.
The only $p$-column of $\bar\kappa$ that is not a $p$-column of $\kappa$ has maximal hook length $m-p$.
Thus there exist sets $A\subseteq H^c(m)$ and $B\subseteq H^c(m-p)$ with $H_{n,p}(\kappa)=(H_{n,p}(\bar\kappa)\cup A)-B$.
Similarly $H_{p,n}(\kappa)=(H_{p,n}(\bar\kappa)\cup C)-D$ for some sets $C\subseteq H^c(m)$ and $D\subseteq H^c(m-n)$.
It suffices to show that $A-B=C-D$ and $B-A=D-C$.

Suppose $z\in A$ but $z\notin B$. Then $z\in H^c(m)$ and $z+n\notin H^c(m)$. On the one hand we obtain $z\notin H^c(m-n)$ and $z\notin D$. It follows that $A\cap D=\emptyset$. On the other hand we obtain $z+n+p\notin H^c(m)$ and therefore $z+n\notin H^c(m-p)$. Since $z\notin B$ this implies $z\notin H^c(m-p)$ and consequently $z+p\notin H^c(m)$. We obtain $z\in C$. Therefore $A-B\subseteq C-D$, and $A-B=C-D$ by symmetry.

Conversely suppose $z\in B$ but $z\notin A$. By symmetry we have $B\cap C=\emptyset$ and $z\notin C$. On the other hand $z+n\notin H^c(m-p)$ implies $z+n+p\notin H^c(m)$ and thus $z+p\notin H^c(m-n)$. Moreover $z\in H^c(m-p)$ implies $z+p\in H^c(m)$ and therefore $z\in H^c(m)$. Since $z\notin A$ we obtain $z+n\in H^c(m)$ and $z\in H^c(m-n)$. We conclude that $z\in D$ and the proof is complete.
\end{proof}

The \emph{skew length} $\op{skl}(\kappa)$ of an $n,p$-core was defined by Armstrong, Hanusa and Jones~\cite{AHJ2014} as the number of cells that are contained in an $n$-row of $\kappa$ and have hook length less than $p$. The set $H_{n,p}(\kappa)$ allows for a new equivalent definition.

\begin{myprop}{skl} Let $\kappa\in\mathfrak C_{n,p}$ be an $n,p$-core. Then $\op{skl}(\kappa)=\#H_{n,p}(\kappa)$.
\end{myprop}

\begin{proof} Fix an $n$-row with largest hook length $h$. On the one hand by \refl{flush} a cell $x$ in this row has hook length less than $p$ if and only if $h_{\kappa}(x)$ is the minimal representative of its residue class modulo $p$ in $H^r(h)$. On the other hand $x$ is contained in a $p$-column if and only if $h_{\kappa}(x)$ is the maximal representative of its residue class modulo $p$ in $H^r(h)$. Thus both $\op{skl}(\kappa)$ and $\#H_{n,p}(\kappa)$ count the number of residue classes modulo $p$ with a representative in $H^r(h)$.
\end{proof}

From \reft{Hab} we immediately recover two results that were recently proven by Guoce Xin~\cite{Xin2015}, and independently by Ceballos, Denton and Hanusa~\cite{CDH2015} when $n$ and $p$ are relatively prime.

\begin{mycor}{skewsym1} The skew length of an $n,p$-core is independent of the order of $n$ and $p$.
\end{mycor}

\begin{mycor}{skewsym2} Let $\kappa\in\mathfrak C_{n,p}$ be an $n,p$-core with conjugate $\kappa'$. Then $\op{skl}(\kappa)=\op{skl}(\kappa')$.
\end{mycor}

Indeed, with our alternative definition of skew length given in \refp{skl} the statements of Corollaries~\ref{Corollary:skewsym1} and~\ref{Corollary:skewsym2} are identical.

\smallskip
We close this section by stating two equivalent conjectural properties of the multiset $H_{n,p}(\kappa)$ that need further investigation.

\begin{myconj}{inH} Let $n$ and $p$ be relatively prime, and $\kappa\in\mathfrak{C}_{n,p}$.
Then the hook length of each cell of $\kappa$ appears in $H_{n,p}(\kappa)$ with multiplicity at least one.
\end{myconj}

\begin{myconj}{invH} Let $n$ and $p$ by relatively prime, $z\geq0$ and $\kappa\in\mathfrak{C}_{n,p}$. Then $z+n\in H_{n,p}(\kappa)$ implies $z\in H_{n,p}(\kappa)$.
\end{myconj}

\section{Rational Shi tableaux}\label{Section:shitab}

In the beginning of this section we recall some facts about root systems and Weyl groups. For further details we refer the reader to~\cite{Humphreys}.
\sk

Let $\Phi$ be an irreducible crystallographic root system with ambient space $V$, positive system $\Phi^+$ and simple system $\Delta=\{\sigma_1,\dots,\sigma_r\}$.
Any root $\alpha\in\Phi$ can be written as a unique integer linear combination
$\alpha=\sum_{i=1}^rc_i\sigma_i$,
where all coefficients $c_i$ are non-negative if $\alpha\in\Phi^+$, or all coefficients are non-positive if $\alpha\in-\Phi^+$.
We define the \emph{height} of the root $\alpha$ by
$\op{ht}(\alpha)=\sum_{i=1}^rc_i$.
Thereby $\op{ht}(\alpha)>0$ if and only if $\alpha\in\Phi^+$, and $\op{ht}(\alpha)=1$ if and only if $\alpha\in\Delta$. Moreover there exists a unique \emph{highest root} $\tilde\alpha$ such that $\op{ht}(\tilde\alpha)\geq\op{ht}(\alpha)$ for all $\alpha\in\Phi$. The \emph{Coxeter number} of $\Phi$ can be defined as $h=\op{ht}(\tilde\alpha)+1$.

\smallskip
Let $\delta$ be a formal variable. We define the set of \emph{affine roots} as
\begin{align*}
\widetilde\Phi
=\{\alpha+k\delta:\alpha\in\Phi,k\in\Z\}\subseteq V\oplus\R\delta.
\end{align*}
The \emph{height} of an affine root is given by $\op{ht}(\alpha+k\delta)=\op{ht}(\alpha)+kh$. The sets of \emph{positive} and \emph{simple affine roots} are defined as
\begin{align*}
\widetilde\Phi^+
=\Phi^+\cup\{\alpha+k\delta\in\widetilde\Phi:\alpha\in\Phi,k>0\}
&&\text{and}&&
\widetilde\Delta
=\Delta\cup\{-\tilde\alpha+\delta\}.
\end{align*}
Thus $\alpha+k\delta\in\widetilde\Phi^+$ if and only if $\op{ht}(\alpha+k\delta)>0$, and $\alpha+k\delta\in\widetilde\Delta$ if and only if $\op{ht}(\alpha+k\delta)=1$.

\smallskip
The \emph{Coxeter arrangement} $\op{Cox}(\Phi)$ consists of all hyperplanes of the form 
$H_{\alpha}=\{x\in V:\skal{x,\alpha}=0\}$
for $\alpha\in\Phi$. Its regions, that is, the connected components of $V-\bigcup_{\alpha\in\Phi}H_{\alpha}$, are called \emph{chambers}. We define the \emph{dominant chamber} as
\begin{align*}
C=\big\{x\in V:\skal{x,\alpha}>0\text{ for all }\alpha\in\Delta\big\}.
\end{align*}
The \emph{Weyl group} $W$ of $\Phi$ is the group of linear automorphisms of $V$ generated by all reflections in a hyperplane in $\op{Cox}(\Phi)$. The Weyl group acts simply transitively on the chambers. Thus identifying the identity $e\in W$ with the dominant chamber, each chamber corresponds to a unique Weyl group element.

The \emph{affine arrangement} $\op{Aff}(\Phi)$ consists of all hyperplanes of the form
$H_{\alpha,k}=\{x\in V:\skal{x,\alpha}=k\}$,
where $\alpha\in\Phi$ and $k\in\Z$.
Its regions are called \emph{alcoves}. We define the \emph{fundamental alcove} as
\begin{align*}
\ac
=\big\{x\in V:\skal{x,\alpha}>0\text{ for all }\alpha\in\Delta\text{ and }\skal{x,\tilde\alpha}<1\big\}.
\end{align*}
The \emph{affine Weyl group} $\W$ of $\Phi$ is the group of affine transformations of $V$ that is generated by all reflections in a hyperplane in $\op{Aff}(\Phi)$.
The affine Weyl group acts simply transitively on the set of alcoves.
By identifying the identity $e\in\W$ with the fundamental alcove, every alcove corresponds to a unique element of $\W$. An element $\omega\in\W$ is called \emph{dominant} if and only if the alcove $\omega(\ac)$ is contained in the dominant chamber $C$.

\smallskip
Given a root $\alpha\in\Phi$ we define its \emph{coroot} as $\alpha^{\vee}=2\alpha/\skal{\alpha,\alpha}$. The \emph{coroot lattice} is the integer span of all coroots
$\coQ=\sum_{\alpha\in\Phi}\Z\alpha^{\vee}\subseteq V$. 
For each $q\in\coQ$ the translation $t_q:V\to V$ defined by $t_q(x)=x+q$ for all $x\in V$ is an element of the affine Weyl group.
Identifying $\coQ$ with its translation group we obtain $\W=W\ltimes\coQ$.
Note that if $\omega\in\W$ is dominant and $\omega=t_qs$, where $q\in\coQ$ and $s\in W$, then in particular $q$ lies in the closure of the dominant chamber.
Thus $\skal{\alpha,q}\geq0$ for each positive root $\alpha\in\Phi^+$.


\sk
The Weyl group is a parabolic subgroup of the affine Weyl group.
Each element $\omega\in\W$ can be assigned a \emph{length} $\ell(\omega)$ indicating the minimal number of elements of a special set of generators needed to express $\omega$.
Each coset $\omega W\in\W/W$ contains a unique representative of minimal length. These representatives are called \emph{Gra{\ss}mannian}. An element $\omega\in\W$ is Gra{\ss}mannian if and only if $\omega^{-1}$ is dominant. 

\smallskip

Fix a root system of type $A_{n-1}$ by declaring the roots, positive roots and simple roots to be
\begin{align*}
\Phi=\{e_i-e_j:i,j\in[n],i\neq j\},
\Phi^+=\{e_i-e_j:i,j\in[n],i<j\},
\Delta=\{e_i-e_{i+1}:i\in[n-1]\}.
\end{align*}
We now describe a combinatorial model for the affine Weyl group of type $A_{n-1}$.
A detailed exposition is found in~\cite[Sec.~8.3]{BjoBre}.
\smallskip

The \emph{affine symmetric group} $\affS_n$ is the group of bijections $\omega:\Z\to\Z$ such that $\omega(i+n)=\omega(i)+n$ for all $i\in\Z$ and $\omega(1)+\dots+\omega(n)=n(n+1)/2$.
Such a bijection $\omega$ is called \emph{affine permutation}.
Each affine permutation is uniquely determined by its \emph{window}
$[\omega(1),\omega(2),\dots,\omega(n)]$.
The group $\affS_n$ has a set of generators called \emph{simple transpositions} given by
$s_i=[1,\dots,i+1,i,\dots,n]$
for $i\in[n-1]$ and
$s_n=[0,2,\dots,n-1,n+1]$.
The affine symmetric group and the set of simple transpositions form a Coxeter system isomorphic to the the affine Weyl group of type $A_{n-1}$.
The symmetric group $\S_n$ can be seen as the subgroup of $\affS_n$ consisting of all affine permutations whose window is a permutation of $[n]$.
Moreover the symmetric group $\S_n$ is the parabolic subgroup of $\affS_n$ generated by $s_i$ for $i\in[n-1]$.
The \emph{length} $\ell(\omega)$ of an affine permutation $\omega\in\affS_n$ is the minimal number $\ell$ of simple transpositions in an expression of the form $\omega=s_{i_1}s_{i_2}\cdots s_{i_{\ell}}$, where $i_j\in[n]$ for all $j\in[\ell]$.
The \emph{Gra{\ss}mannian affine permutations}, that is, the minimal length representatives of the cosets in $\affS_n/\S_n$, are described in \cite[Prop.~8.3.4]{BjoBre}.

\begin{myprop}{grass} An affine permutation $\omega\in\affS_n$ is the minimal length representative of its coset $\omega\S_n\in\affS_n/\S_n$ if and only if its window is increasing, that is,
\begin{align*}
\omega(1)<\omega(2)<\dots<\omega(n).
\end{align*}
\end{myprop}


We now give a combinatorial description of the decomposition of the affine symmetric group as the semidirect product of translations and permutations.
The coroot lattice of type $A_{n-1}$ is given by
\begin{align*}
\coQ=\big\{x\in\Z^n:x_1+\dots+x_n=0\big\}.
\end{align*}
The affine symmetric group acts on the coroot lattice via the following rules
\begin{align*}
s_i\act x=(x_1,\dots,x_{i+1},x_i,\dots,x_n)&&\text{for }i\in[n-1]\text{ and}&&
s_n\act x=(x_n+1,x_2,\dots,x_{n-1},x_1-1).
\end{align*}
Given $\omega\in\affS_n$ write $\omega(i)=a_in+b_i$ for each $i\in\Z$, where $a_i\in\Z$ and $b_i\in[n]$. Set $\sigma(\omega,i)=b_i$, $\mu(\omega,b_i)=a_i$ and $\nu(\omega,i)=-a_i$.

\begin{mylem}{sn}
\begin{enumerate}[(i)]
\item\label{item:sigma} The assignment $i\mapsto\sigma(\omega,i)$ defines a permutation $\sigma(\omega)\in\S_n$.
\item\label{item:munu} The vectors $\mu(\omega)=(\mu(\omega,1),\dots,\mu(\omega,n))$ and $\nu(\omega)=(\nu(\omega,1),\dots,\nu(\omega,n))$ lie in the coroot lattice $\coQ$ and for all $i\in[n]$ we have
\begin{align*}
\mu(s_i\omega)=s_i\act\mu(\omega),&&
\nu(\omega s_i)=s_i\act\nu(\omega).
\end{align*}
\item\label{item:atzero} We have $\omega\act0=\mu(\omega)$ and $\omega^{-1}\act0=\nu(\omega)$.
\item\label{item:munuinv} We have $\mu(\omega^{-1})=\nu(\omega)$ and $\sigma(\omega^{-1})=\sigma(\omega)^{-1}$.
\item\label{item:munusigma} We have $\mu(\omega)=-\sigma(\omega)\act\nu(\omega)$.
\end{enumerate}
\end{mylem}

\begin{proof} Claims~\refi{sigma} and~\refi{munusigma} are straight forward. Claims~\refi{atzero} and~\refi{munuinv} follow immediately from~\refi{munu}, which can be shown by induction on the length of $\omega$.

Clearly $\mu(e)=\nu(e)=(0,\dots,0)\in\coQ$. Set $\sigma=\sigma(\omega)$ and suppose that $\mu(\omega)\in\coQ$.

Fix $i\in[n-1]$ and choose $j,k\in[n]$ such that $\sigma(j)=i$ and $\sigma(k)=i+1$.
Then $s_i\omega(\ell)=\omega(\ell)$ for all $\ell\in[n]-\{j,k\}$ and hence $\mu(s_i\omega,\ell)=\mu(\omega,\ell)$ for all $\ell\in[n]-\{i,i+1\}$.
Furthermore, $s_i\omega(j)=s_i(a_jn+i)=a_jn+i+1$ and $s_i\omega(k)=s_i(a_kn+i+1)=a_kn+i$.
It follows that $\mu(s_i\omega,i)=a_k=\mu(\omega,i+1)$ and $\mu(s_i\omega,i+1)=a_j=\mu(\omega,i)$.

Next choose $j,k\in[n]$ such that $\sigma(j)=1$ and $\sigma(k)=n$.
Then $s_n\omega(\ell)=\omega(\ell)$ for all $\ell\in[n]-\{j,k\}$.
Thus $\mu(s_n\omega,\ell)=\mu(\omega,\ell)$ for all $\ell\in[n]-\{1,n\}$.
Moreover, $s_n\omega(j)=s_n(a_jn+1)=a_jn+0=(a_j-1)n+n$ and $s_n\omega(k)=s_n(a_kn+n)=a_kn+(n+1)=(a_k+1)n+1$.
Thus $\mu(s_n\omega,1)=a_k+1=\mu(\omega,n)+1$ and $\mu(s_n\omega,n)=a_j-1=\mu(\omega,1)-1$.

We obtain $\mu(s_i\omega)=s_i\act\mu(\omega)$, and in particular $\mu(s_i\omega)\in\coQ$, for all $i\in[n]$.

Now suppose $\nu(\omega)\in\coQ$ and fix $i\in[n-1]$.
Since $\omega s_i(j)=\omega(j)$ for all $j\in[n]-\{i,i+1\}$ we have $\nu(\omega s_i,j)=\nu(\omega,j)$ for all $j\in[n]-\{i,i+1\}$.
Furthermore $\omega s_i(i)=\omega(i+1)$ and $\omega s_i(i+1)=\omega(i)$, thus $\nu(\omega s_i,i)=\nu(\omega,i+1)$ and $\nu(\omega s_i,i+1)=\nu(\omega,i)$.

Finally since $\omega s_n(j)=\omega(j)$ for all $j\in[n]-\{1,n\}$ we have $\nu(\omega s_n,j)=\nu(\omega,j)$ for all $j\in[n]-\{1,n\}$.
Moreover $\omega s_n(1)=\omega(0)=\omega(n-n)=\omega(n)-n=(a_n-1)n+b_n$, hence $\nu(\omega s_n,1)=-a_n+1=\nu(\omega,n)+1$.
On the other hand $\omega s_n(n)=\omega(n+1)=n+\omega(1)=(a_1+1)n+b_1$.
Thus $\nu(\omega s_n,n)=-a_1-1=\nu(\omega,1)-1$.

We obtain that $\nu(\omega s_i)=s_i\act\nu(\omega)$, and in particular $\nu(\omega s_i)\in\coQ$, for all $i\in[n]$. This concludes the proof of claim \refi{munu} and thus the proof of the lemma.
\end{proof}

For $q\in\coQ$ define an affine permutation $t_q\in\affS_n$ by $t_q(i)=q_in+i$ for $i\in[n]$. 
We call an affine permutation $\omega\in\affS_n$ a \emph{translation} if there exists $q\in\coQ$ such that $\omega\act x=q+x$ for all $x\in\coQ$.

\begin{mythrm}{translation}
\begin{enumerate}[(i)]
\item\label{item:tmusigma} Let $\omega\in\affS_n$ be an affine permutation and set $s=\sigma(\omega)$, $x=\mu(\omega)$ and $y=\nu(\omega)$. Then
$\omega=t_xs=st_{-y}$.
\item\label{item:coQ} Let $x,y\in\coQ$. Then $t_xt_y=t_{x+y}$ and $(t_x)^{-1}=t_{-x}$. Hence we may view the coroot lattice $\coQ$ as a subgroup of $\affS_n$.
\item\label{item:translation} An affine permutation $\omega\in\affS_n$ is a translation if and only if $\omega=t_q$ for some $q\in\coQ$
\item\label{item:semidirect} The affine symmetric group is the semidirect product of the symmetric group and the coroot lattice, that is,
$\affS_n=\S_n\ltimes\coQ$.
\end{enumerate}
\end{mythrm}

\begin{proof} Claims~\refi{tmusigma} and~\refi{coQ} are straightforward calculations.
Claim~\refi{tmusigma} follows from
\begin{align*}
t_xs(i)
&=\mu(\omega,s(i))n+s(i)
=\omega(i)
=-\nu(\omega,i)n+s(i)
=s(-\nu(\omega,i)n+i)
=st_{-y}(i).
\end{align*}
On the other hand
\begin{align*}
t_xt_y(i)
&=t_x(y_in+i)
=y_in+t_x(i)
=(x_i+y_i)n+i
=t_{x+y}(i)
\end{align*}
implies~\refi{coQ}.
To see~\refi{translation} note that $t_q\act x=t_qt_x\act0=t_{q+x}\act0=q+x$. 
Conversely, if $\omega$ is a translation by $q\in\coQ$ then $q+x=\omega\act x=q+\sigma(\omega)\act x$ for all $x\in\coQ$, which implies $\sigma(\omega)=e$. 
Finally, $\S_n\coQ=\affS_n$ and $\S_n$ normalises $\coQ$ by~\refi{tmusigma}. Since $\S_n\cap \coQ=\{e\}$, we obtain~\refi{semidirect}.
\end{proof}

Thus the decomposition of an affine permutation into a product of a translation and a permutation can be obtained from its window in a simple and direct fashion.
We remark that this combinatorial decomposition has appeared in the literature before, for example in~\cite{BjoBre1996}. However, the author is unaware of a reference explaining explicitly its connection to the algebraic decomposition into a semidirect product, which exists for any affine Weyl group.

\smallskip
The affine symmetric group possesses an involutive automorphism owing to the symmetry of the Dynkin diagram of type $A_{n-1}$. Set $s_i^*=s_{n-i}$ for $i\in[n-1]$ and $s_n^*=s_n$. This correspondence extends to an automorphism $\omega\mapsto\omega^*$ on $\affS_n$, where $\omega^*$ is obtained by replacing all instances of $s_i$ in any expression of $\omega$ in terms of the simple transpositions by $s_i^*$. The involutive automorphism has a simple explicit description in window notation and fulfils many desirable properties.
Some of these properties are presented with proofs in this paper, although they are already well-known to experts.

\begin{mylem}{aut}
\begin{enumerate}[(i)]
\item\label{item:window} Let $\omega\in\affS_n$ be an affine permutation and $i\in\Z$. Then $w^*(i)=1-\omega(1-i)$. In particular the window of $\omega^*$ is given by
$[n+1-\omega(n),\dots,n+1-\omega(1)]$.
\item\label{item:aut} The involutive automorphism preserves $\S_n$, translations, Gra{\ss}mannian affine permutations and dominant affine permutations.
\end{enumerate}
\end{mylem}

\begin{proof} We prove claim~\refi{window} by induction on the length of $\omega$. Clearly~\refi{window} holds for the identity, as $e^*=e$.
Thus assume it is true for $\omega$.
For $i\in[n-1]$ the right multiplication of $\omega$ by $s_i$ corresponds to exchanging the two numbers $\omega(i),\omega(i+1)$ in the window of $\omega$.
On the other hand multiplying $\omega^*$ by $s_i^*$ from the right exchanges the numbers $n+1-\omega(i+1),n+1-\omega(i)$. Claim~\refi{window} therefore also holds for $\omega s_i$.

Claim~\refi{aut} follows from~\refi{window} and the fact that $(\omega^{-1})^*=(\omega^*)^{-1}$.
\end{proof}

%

\begin{figure}[t]
\begin{center}
\begin{tikzpicture}[scale=.6]
\begin{scope}
\draw[xshift=5mm,yshift=5mm]
(0,0) node{\footnotesize{$k_{1,2}$}}
(0,1) node{\footnotesize{$k_{1,3}$}}
(0,2) node{\footnotesize{$k_{1,4}$}}
(0,3) node{\footnotesize{$k_{1,5}$}}
(0,4) node{\footnotesize{$k_{1,6}$}}
(0,5) node{\footnotesize{$k_{1,7}$}}
(1,1) node{\footnotesize{$k_{2,3}$}}
(1,2) node{\footnotesize{$k_{2,4}$}}
(1,3) node{\footnotesize{$k_{2,5}$}}
(1,4) node{\footnotesize{$k_{2,6}$}}
(1,5) node{\footnotesize{$k_{2,7}$}}
(2,2) node{\footnotesize{$k_{3,4}$}}
(2,3) node{\footnotesize{$k_{3,5}$}}
(2,4) node{\footnotesize{$k_{3,6}$}}
(2,5) node{\footnotesize{$k_{3,7}$}}
(3,3) node{\footnotesize{$k_{4,5}$}}
(3,4) node{\footnotesize{$k_{4,6}$}}
(3,5) node{\footnotesize{$k_{4,7}$}}
(4,4) node{\footnotesize{$k_{5,6}$}}
(4,5) node{\footnotesize{$k_{5,7}$}}
(5,5) node{\footnotesize{$k_{6,7}$}}
;
\end{scope}
\draw(7,3)node{$=$};
\begin{scope}[xshift=8cm]
\draw[xshift=5mm,yshift=5mm]
(0,0) node{\footnotesize{$0$}}
(0,1) node{\footnotesize{$1$}}
(0,2) node{\footnotesize{$2$}}
(0,3) node{\footnotesize{$2$}}
(0,4) node{\footnotesize{$3$}}
(0,5) node{\footnotesize{$5$}}
(1,1) node{\footnotesize{$1$}}
(1,2) node{\footnotesize{$2$}}
(1,3) node{\footnotesize{$2$}}
(1,4) node{\footnotesize{$2$}}
(1,5) node{\footnotesize{$5$}}
(2,2) node{\footnotesize{$1$}}
(2,3) node{\footnotesize{$1$}}
(2,4) node{\footnotesize{$1$}}
(2,5) node{\footnotesize{$3$}}
(3,3) node{\footnotesize{$0$}}
(3,4) node{\footnotesize{$0$}}
(3,5) node{\footnotesize{$2$}}
(4,4) node{\footnotesize{$0$}}
(4,5) node{\footnotesize{$2$}}
(5,5) node{\footnotesize{$2$}}
;
\end{scope}
\end{tikzpicture}
\caption{The inversion table of the affine permutation $\omega=[-2,15,-1,16,-14,10,4]\in\smash{\stackrel{\sim}{\smash{\mathfrak{S}}\protect\rule{0pt}{1.05ex}}}_7$ with inverse $\omega^{-1}=[-12,-10,-1,7,8,10,26]$.}
\label{Figure:kij}
\end{center}
\end{figure}

Let $\omega\in\affS_n$ be an affine permutation. An \emph{affine inversion} of $\omega$ is a pair $(i,j)\in[n]\times\N$ such that $i<j$ and $\omega(i)>\omega(j)$. For $i,j\in[n]$ with $i<j$ define
\begin{align}\label{eq:kij}
k_{i,j}(\omega)
&=\abs{\floor{\frac{\omega^{-1}(j)-\omega^{-1}(i)}{n}}}.
\end{align}
The numbers $k_{i,j}(\omega)$ where first considered by Shi who proved that $\sum_{i,j}k_{i,j}(\omega)=\ell(\omega)$ for all $\omega\in\affS_n$~\cite[Lem.~4.2.2]{Shi:Kazhdan_Lusztig}.
We arrange these numbers in a staircase tableau as in \reff{kij}.
The following lemma explains that $k_{i,j}(\omega)$ counts certain affine inversions of $\omega$. It should be compared to the discussion preceding~\cite[Prop.~8.3.1]{BjoBre}.
\begin{mylem}{invA} Let $\omega\in\affS_n$ be an affine permutation.
\begin{enumerate}[(i)]
\item\label{item:inv} Then the number of affine inversions of $\omega$ is given by $\sum_{i,j}k_{i,j}(\omega)$. 
\item\label{item:invij} If $\omega$ is dominant then $k_{i,j}(\omega)$ equals the number of affine inversions $(a,b)$ of $\omega$, such that $\omega(a)\equiv i$ and $\omega(b)\equiv j$ modulo $n$.
\end{enumerate}
\end{mylem}

\begin{proof} Let $i,j\in[n]$ with $i<j$. If $\omega^{-1}(i)<\omega^{-1}(j)$ then $(j,i+kn)$ is an affine inversion of $\omega^{-1}$ if and only if $k\geq1$ and
\begin{align*}
\omega^{-1}(j)>\omega^{-1}(i+kn)
&\Leftrightarrow\omega^{-1}(j)-\omega^{-1}(i)>kn
\Leftrightarrow\floor{\frac{\omega^{-1}(j)-\omega^{-1}(i)}{n}}\geq k.
\end{align*}
On the other and if $\omega^{-1}(i)>\omega^{-1}(j)$ then $(i,j+kn)$ is an affine inversion of $\omega^{-1}$ if and only if $k\geq0$ and
\begin{align*}
k\leq\floor{\frac{\omega^{-1}(i)-\omega^{-1}(j)}{n}}
=\abs{\floor{\frac{\omega^{-1}(j)-\omega^{-1}(i)}{n}}}-1.
\end{align*}
This proves the first claim since $\omega$ and $\omega^{-1}$ have the same number of inversions.

If $\omega$ is dominant then $\omega^{-1}$ is Gra{\ss}mannian. Thus $\omega^{-1}(i)<\omega^{-1}(j)$ and $k_{i,j}(\omega)$ counts the number of inversions of $\omega^{-1}$ of the form $(j,i+kn)$ as above. But these inversions correspond exactly to the inversions $(a,b)$ of $\omega$ with $\omega(a)\equiv i$ and $\omega(b)\equiv j$ modulo $n$. Claim~\refi{invij} follows.
\end{proof}

For example the dominant affine permutation $\omega=[-2,15,-1,16,-14,10,4]$ from \reff{kij} has three affine inversions $(a,b)$ such that $\omega(a)\equiv3$ and $\omega(b)\equiv7$ modulo $7$, namely $(6,12),(6,19)$ and $(6,26)$. Thus $k_{3,7}(\omega)=3$.

Motivated by \refl{invA} we call the collection of numbers $k_{i,j}(\omega)$ for $i,j\in[n]$ with $i<j$ the \emph{inversion table} of $\omega\in\affS_n$.
A simple computation reveals the effect of the involutive automorphism on the inversion table.

\begin{myprop}{invconj} Let $\omega\in\affS_n$ be an affine permutation. Then the inversion table of $\omega^*$ is the transpose of the inversion table of $\omega$.
That is,
$k_{n+1-j,n+1-i}(\omega^*)=k_{i,j}(\omega)$
for all $i,j\in[n]$ with $i<j$.
\end{myprop}

\begin{proof} We compute
\begin{align*}
k_{n+1-j,n+1-i}(\omega^*)
&=\abs{\floor{\frac{(\omega^*)^{-1}(n+1-i)-(\omega^*)^{-1}(n+1-j)}{n}}}\\
&=\abs{\floor{\frac{(\omega^{-1})^*(n+1-i)-(\omega^{-1})^*(n+1-j)}{n}}}\\
&=\abs{\floor{\frac{n+1-\omega^{-1}(i)-(n+1)+\omega^{-1}(j)}{n}}}
=k_{i,j}(\omega).\qedhere
\end{align*}
\end{proof}

There is a natural way to generalise the inversion table to affine Weyl groups. Let $\Phi$ be an irreducible crystallographic root system with affine Weyl group $\W$. Let $\omega\in\W$ have the decomposition $\omega=t_qs$, where $q\in\coQ$ and $s\in W$. The affine Weyl group acts on the set of affine roots by
\begin{align*}
\omega\act(\alpha+k\delta)
=s\act\alpha+(k-\skal{q,s\act\alpha})\delta.
\end{align*}
A positive affine root $\alpha+k\delta\in\widetilde{\Phi}^+$ is an \emph{affine inversion} of $\omega$ if $\omega\act(\alpha+k\delta)\in-\widetilde{\Phi}^+$.
Denote the set of affine inversions of $\omega$ by
$\op{Inv}(\omega)=\widetilde{\Phi}^+\cap\omega^{-1}\act\big(-\widetilde{\Phi}^+\big)$.
We remark that the affine inversions of $\omega$ correspond to the affine hyperplanes in $\op{Aff}(\Phi)$ separating $\omega(\ac)$ from the fundamental alcove.
Thus the number of inversions equals the length of $\omega$.

For $\omega\in\W$ and a positive root $\alpha\in\Phi^+$ define
\begin{align}\label{eq:kalpha}
k_{\alpha}(\omega)
&=\#(\op{Inv}(\omega^{-1})\cap\{\pm\alpha+k\delta:k\in\Z\}).
\end{align}

\begin{mylem}{kalpha} Let $\alpha\in\Phi^+$ be a positive root and $\omega\in\W$ an element of the affine Weyl group. If $\omega=t_qs$, where $q\in\coQ$ and $s\in W$, then
\begin{align*}
k_{\alpha}(\omega)=
\begin{cases}
\abs{\skal{\alpha,q}}&\quad\text{if }s^{-1}\act\alpha\in\Phi^+,\\
\abs{\skal{\alpha,q}-1}&\quad\text{if }s^{-1}\act\alpha\in-\Phi^+.
\end{cases}
\end{align*}
\end{mylem}

\begin{proof} Set $\beta=s^{-1}\act\alpha$. For $k\geq0$
\begin{align*}
\omega^{-1}\act(\alpha+k\delta)
&=\beta+(k+\skal{\alpha,q})\delta
\in-\widetilde\Phi^+
\end{align*}
if and only if either $k<-\skal{\alpha,q}$ or both $k=-\skal{\alpha,q}$ and $\beta\in-\Phi^+$. On the other hand for $k\geq1$
\begin{align*}
\omega^{-1}\act(-\alpha+k\delta)
&=-\beta+(k-\skal{\alpha,q})\delta
\in-\widetilde\Phi^+
\end{align*}
if and only if $k<\skal{\alpha,q}$ or both $k=\skal{\alpha,q}$ and $\beta\in\Phi^+$. Combined this implies the claim.
\end{proof}

As a consequence we deduce that the collection of numbers $k_{\alpha}(\omega)$ defined in \refq{kalpha} really is a generalisation of the inversion table defined in \refq{kij} for type $A_{n-1}$.

\begin{myprop}{kalpha} Let $i,j\in[n]$ with $i<j$ and $\omega\in\affS_n$. Then $k_{i,j}(\omega)=k_{e_i-e_j}(\omega)$.
\end{myprop}

\begin{proof} Write $\omega=t_qs$, where $q\in\coQ$ and $s\in\S_n$. The claim follows from \refl{kalpha} and
\begin{align*}
\abs{\floor{\frac{\omega^{-1}(j)-\omega^{-1}(i)}{n}}}
&=\abs{\floor{\frac{-q_jn+s^{-1}(j)-(-q_in+s^{-1}(i))}{n}}}\\
&=
\begin{cases}
\abs{q_i-q_j}&\quad\text{if }s^{-1}(i)<s^{-1}(j),\\
\abs{q_i-q_j-1}&\quad\text{if }s^{-1}(i)>s^{-1}(j).
\end{cases}\qedhere
\end{align*}

\end{proof}

Note that two distinct elements of $\W$ may well have the same inversion table. For example the two affine permutations $[0,1,5]$ and $[-1,3,4]$ in $\affS_3$ have the same inversion table. However, the inversion table is closely related to another construction studied by Shi. Let $\omega\in\W$ and $\alpha\in\Phi^+$. Define $\tilde{k}_{\alpha}(\omega)\in\Z$ such that 
\begin{align*}
\tilde{k}_{\alpha}(\omega)
<\skal{\alpha,x}
<\tilde{k}_{\alpha}(\omega)+1
\end{align*}
for all $x\in\omega(\ac)$. The collection of numbers $\tilde{k}_{\alpha}(\omega)$ for $\alpha\in\Phi^+$ is called \emph{address} or \emph{Shi coordinates} of the alcove $\omega(\ac)$. Clearly, each alcove is determined uniquely by its Shi coordinates, and $\tilde{k}_{\alpha}(\omega)\geq0$ for all $\alpha\in\Phi^+$ if and only if $\omega$ is dominant. Moreover it is not difficult to show that
\begin{align*}
\tilde{k}_{\alpha}(\omega)=
\begin{cases}
\skal{\alpha,q}&\quad\text{if }s^{-1}\act\alpha\in\Phi^+,\\
\skal{\alpha,q}-1&\quad\text{if }s^{-1}\act\alpha\in-\Phi^+,
\end{cases}
\end{align*}
where $q\in\coQ$ and $s\in W$ are chosen such that $\omega=t_qs$ (See~\cite[Thm.~3.3]{Shi1987:alcoves}).
%
%
%
Thus we obtain $k_{\alpha}(\omega)=\abs{\smash{\tilde{k}}_{\alpha}(\omega)}$ for all $\omega\in\W$ and $\alpha\in\Phi^+$.
In particular it follows that no two distinct dominant\,(!) elements of the affine Weyl group have the same inversion table.

\sk
Similar to the inversion table, the entries of the rational Shi tableau count certain inversions of an affine permutation.
However, as we have just seen, the set of all affine inversions might be too big to obtain a nice correspondence between permutations and tableaux.
In case of the inversion table it is more fruitful to restrict ourselves to dominant affine permutations.
In case of the rational Shi tableau the appropriate set is even smaller, in fact finite.

Let $p$ be a positive integer that is relatively prime to $n$.
An affine permutation $\omega\in\affS_n$ is called \emph{$p$-stable} if $\omega(i)<\omega(i+p)$ for all $i\in\Z$.
The set of $p$-stable affine permutations in $\affS_n$, denoted by $\affS_n^p$, was first considered by Gorsky, Mazin and Vazirani~\cite{GMV2014}.
As a consequence of~\refl{aut} we see that the involutive automorphism $\omega\mapsto\omega^*$ preserves $\affS_n^p$.

If $p=mn+1$ then $\affS_n^p$ corresponds to the set of minimal alcoves of the regions of the $m$-extended Shi arrangement. If $p=mn-1$ then $\affS_n^p$ is related to the bounded regions of the $m$-extended Shi arrangement. Thus $p$-stable affine permutations should be viewed as a ``rational'' analogue of the regions of the Shi arrangement.

Let $\omega\in\affS_n^p$ be a dominant $p$-stable affine permutation.
Let $\omega^{-1}=t_qs$, where $q\in\coQ$ and $s\in\S_n$, and $p=mn+r$ with $r\in[n-1]$.
We define the \emph{rational Shi tableau} of $\omega$ as the collection of integers
\begin{align}\label{eq:tij}
t^p_{i,j}(\omega)=
\begin{cases}
\min(k_{i,j}(\omega),m)&\quad\text{if }r+s(i)<s(j)\text{ or }s(i)+r-n<s(j)<s(i).\\
\min(k_{i,j}(\omega),m+1)&\quad\text{otherwise,}
\end{cases}
\end{align}
where $i,j\in[n]$ such that $i<j$. See \reff{tij} for an example.

\smallskip
Note that $k_{i,j}(\omega)\leq m$ whenever $s(j)-s(i)\in\{r-n,r\}$.
This follows from \refl{kalpha} and the fact that $\omega$ is dominant $p$-stable.
As a consequence $t_{i,j}^p(\omega)=\min(k_{i,j}(\omega),m)$ if $p=mn+1$, and $t_{i,j}^p(\omega)=\min(k_{i,j}(\omega),m+1)$ if $p=mn+(n-1)$.
For $p=mn+1$ we recover the Shi tableau of Fishel, Tzanaki and Vazirani~\cite{FTV2011}

\begin{figure}[t]
\begin{center}
\begin{tikzpicture}[scale=.6]
\begin{scope}
\draw[xshift=5mm,yshift=5mm]
(0,0) node{\footnotesize{$t_{1,2}$}}
(0,1) node{\footnotesize{$t_{1,3}$}}
(0,2) node{\footnotesize{$t_{1,4}$}}
(0,3) node{\footnotesize{$t_{1,5}$}}
(0,4) node{\footnotesize{$t_{1,6}$}}
(0,5) node{\footnotesize{$t_{1,7}$}}
(1,1) node{\footnotesize{$t_{2,3}$}}
(1,2) node{\footnotesize{$t_{2,4}$}}
(1,3) node{\footnotesize{$t_{2,5}$}}
(1,4) node{\footnotesize{$t_{2,6}$}}
(1,5) node{\footnotesize{$t_{2,7}$}}
(2,2) node{\footnotesize{$t_{3,4}$}}
(2,3) node{\footnotesize{$t_{3,5}$}}
(2,4) node{\footnotesize{$t_{3,6}$}}
(2,5) node{\footnotesize{$t_{3,7}$}}
(3,3) node{\footnotesize{$t_{4,5}$}}
(3,4) node{\footnotesize{$t_{4,6}$}}
(3,5) node{\footnotesize{$t_{4,7}$}}
(4,4) node{\footnotesize{$t_{5,6}$}}
(4,5) node{\footnotesize{$t_{5,7}$}}
(5,5) node{\footnotesize{$t_{6,7}$}}
;
\end{scope}
\draw(7,3)node{$=$};
\begin{scope}[xshift=8cm]
\draw[xshift=5mm,yshift=5mm]
(0,0) node[shape=circle,draw,inner sep=.1pt,minimum size=1.2em]{\footnotesize{$0$}}
(0,1) node{\footnotesize{$1$}}
(0,2) node{\footnotesize{$2$}}
(0,3) node{\footnotesize{$2$}}
(0,4) node[shape=circle,draw,inner sep=.1pt,minimum size=1.2em]{\footnotesize{$3$}}
(0,5) node{\footnotesize{$2$}}
(1,1) node[shape=circle,draw,inner sep=.1pt,minimum size=1.2em]{\footnotesize{$1$}}
(1,2) node{\footnotesize{$2$}}
(1,3) node{\footnotesize{$2$}}
(1,4) node{\footnotesize{$2$}}
(1,5) node[shape=circle,draw,inner sep=.1pt,minimum size=1.2em]{\footnotesize{$3$}}
(2,2) node[shape=circle,draw,inner sep=.1pt,minimum size=1.2em]{\footnotesize{$1$}}
(2,3) node[shape=circle,draw,inner sep=.1pt,minimum size=1.2em]{\footnotesize{$1$}}
(2,4) node{\footnotesize{$1$}}
(2,5) node{\footnotesize{$2$}}
(3,3) node[shape=circle,draw,inner sep=.1pt,minimum size=1.2em]{\footnotesize{$0$}}
(3,4) node{\footnotesize{$0$}}
(3,5) node{\footnotesize{$2$}}
(4,4) node[shape=circle,draw,inner sep=.1pt,minimum size=1.2em]{\footnotesize{$0$}}
(4,5) node{\footnotesize{$2$}}
(5,5) node[shape=circle,draw,inner sep=.1pt,minimum size=1.2em]{\footnotesize{$2$}}
;
\end{scope}
\end{tikzpicture}
\caption{The rational Shi tableau of $\omega=[-2,15,-1,16,-14,10,4]\in\smash{\stackrel{\sim}{\smash{\mathfrak{S}}\protect\rule{0pt}{1.05ex}}}_7^{16}$.
We have $m=2$, $r=2$ and $s=[5,1,6,2,7,3,4]$.
The circled entries indicate that a minimum with $m+1=3$ is taken.
Also compare to \reff{kij}.}
\label{Figure:tij}
\end{center}
\end{figure}

\sk
The following proposition shows that the rational Shi tableau behaves similarly to the inversion table under the involutive automorphism.

\begin{myprop}{shiconj} Let $\omega\in\affS_n^p$ be a dominant $p$-stable affine permutation. Then the rational Shi tableau of $\omega^*$ is the transpose of the rational Shi tableau of $\omega$. That is,
$t^p_{n+1-j,n+1-i}(\omega^*)=t^p_{i,j}(\omega)$
for all $i,j\in[n]$ with $i<j$.
\end{myprop}

\begin{proof} The claim follows from \refp{invconj} and
\begin{align*}
r+s(i)<s(j)
&\Leftrightarrow r+n+1-s^*(n+1-i)<n+1-s^*(n+1-j)\\
&\Leftrightarrow r+s^*(n+1-j)<s^*(n+1-i),\\
s(i)+r-n<s(j)
&\Leftrightarrow n+1-s^*(n+1-i)+r-n<n+1-s^*(n+1-j)\\
&\Leftrightarrow s^*(n+1-j)+r-n<s^*(n+1-i),\\
s(j)<s(i)
&\Leftrightarrow s^*(n+1-i)<s^*(n+1-j).\qedhere
\end{align*}
\end{proof}

The rational Shi tableau generalises naturally to affine Weyl groups.
Let $\Phi$ be an irreducible crystallographic root system, and $p$ be relatively prime to the Coxeter number $h$ of $\Phi$.
Let $\widetilde\Phi_p$ denote the set of affine roots of height $p$.
In~\cite[Sec.~8]{Thiel2015} Thiel defines \emph{$p$-stable} affine Weyl group elements as the set 
\begin{align*}
\W^p=\big\{\omega\in\W:\omega\act\widetilde\Phi_p\subseteq\widetilde\Phi^+\big\}.
\end{align*}
Let $\alpha\in\Phi^+$ be a positive root and $\omega\in\W^p$ a dominant $p$-stable element of the affine Weyl group. \smallskip
Define
\begin{align}\label{eq:talpha}
t_{\alpha}^p(\omega)
&=\#\big(\omega\act(-\widetilde\Phi_{<p}^+)\cap\{-\alpha+k\delta:k\geq1\}\big),
\end{align}
where $\widetilde\Phi^+_{<p}$ denotes the set of positive affine roots with height less than $p$.
The \emph{rational Shi tableau} of $\omega$ is the collection of numbers $t_{\alpha}^p(\omega)$ for $\alpha\in\Phi^+$.

\begin{mylem}{dominant} Let $\alpha\in\Phi^+$ be a positive root and $\omega\in\W$ be dominant with $\omega=t_qs$, where $q\in\coQ$ and $s\in W$. If $\skal{\alpha,q}=0$ then $s^{-1}\act\alpha\in\Phi^+$.
\end{mylem}

\begin{proof} Suppose $\skal{\alpha,q}=0$. The height function $\op{ht}:\Phi\to\R$ extends to a linear functional on $V$. Thus we may chose $v\in V$ with $\skal{\beta,v}=\op{ht}(\beta)/h$ for all $\beta\in\Phi$, where $h$ is the Coxeter number of $\Phi$. Note that $v\in\ac$ by definition. Thus $\skal{\alpha,\omega\act v}>0$ since $\omega$ is dominant. We compute
\begin{align*}
\frac{\op{ht}(s^{-1}\act\alpha)}{h}
&=\skal{s^{-1}\act\alpha,v}
=\skal{\alpha,s\act v}
=\skal{\alpha,q+s\act v}
=\skal{\alpha,\omega\act v}
>0.\qedhere
\end{align*}
\end{proof}

Using \refl{dominant} we obtain that $\omega^{-1}\act(\alpha+k\delta)\in\widetilde\Phi^+$ for all $\alpha\in\Phi^+$ and $k\geq0$ whenever $\omega$ is dominant. Hence,
\begin{align*}
t_{\alpha}^p(\omega)
=\#(\op{Inv}(\omega^{-1})\cap\{\pm\alpha+k\delta:k\in\Z\}\cap\omega\act(-\widetilde\Phi_{<p})).
\end{align*}
We see that similar to $k_{\alpha}(\omega)$ also $t_{\alpha}^p(\omega)$ counts certain affine inversions $\pm\alpha+k\delta$ of $\omega^{-1}$, but with an additional restriction on the height of $\omega^{-1}\act(\pm\alpha+k\delta)$.

\begin{myprop}{talpha} Let $m\geq0$ and $r\in[h-1]$, where $h$ is the Coxeter number of $\Phi$, and $\alpha\in\Phi^+$ be a positive root. Set $p=mh+r$ and let $\omega\in\W^p$ be dominant $p$-stable. If $\omega=t_qs$, where $q\in\coQ$ and $s\in W$, then
\begin{align*}
t_{\alpha}^p(\omega)
&=
\begin{cases}
\min(k_{\alpha}(\omega),m)&\quad\text{if }r-h<\op{ht}(s^{-1}\act\alpha)<0\text{ or }r<\op{ht}(s^{-1}\act\alpha),\\
\min(k_{\alpha}(\omega),m+1)&\quad\text{otherwise.}
\end{cases}
\end{align*}
\end{myprop}

\begin{proof} Set $\beta=s^{-1}\act\alpha$. Recall that $\skal{\alpha,q}\geq0$ because $\omega$ is dominant. 
If $\op{ht}(\beta)=r$ then $\beta+m\delta\in\widetilde\Phi_p$ and therefore
\begin{align*}
\alpha+(m-\skal{\alpha,q})\delta
&=\omega\act(\beta+m\delta)
\in\widetilde\Phi^+
\end{align*}
since $\omega\in\W^p$. If instead $\op{ht}(\beta)=r-h$ then $\beta+(m+1)\delta\in\widetilde\Phi_p$ and
\begin{align*}
\alpha+(m+1-\skal{\alpha,q})\delta
&=\omega\act(\beta+(m+1)\delta)
\in\widetilde\Phi^+.
\end{align*}
Consequently $\op{ht}(\beta)=r$ implies $\skal{\alpha,q}\leq m$, and $\op{ht}(\beta)=r-h$ implies $\skal{\alpha,q}\leq m+1$.
Now let $k\geq1$.
Then
\begin{align*}
0<\op{ht}\big(-\omega^{-1}\act(-\alpha+k\delta)\big)
=\op{ht}(\beta)+(\skal{\alpha,q}-k)h<p=mh+r
\end{align*}
if and only if one of the following (mutually exclusive) cases occurs
\begin{align*}
m>0\text{ and }k&=\skal{\alpha,q}\text{ and }0<\op{ht}(\beta),\\
m=0\text{ and }k&=\skal{\alpha,q}\text{ and }0<\op{ht}(\beta)<r,\\
-m&<k-\skal{\alpha,q}<0,\\
m>0\text{ and }-m&=k-\skal{\alpha,q}\text{ and }\op{ht}(\beta)<r,\\
-m-1&=k-\skal{\alpha,q}\text{ and }\op{ht}(\beta)<-h+r.
\end{align*}
Equivalently $-\alpha+k\delta$ contributes to $t_{\alpha}^p(\omega)$ if and only if
\begin{align*}
k\in\{1,2,\dots\}\cap
\begin{cases}
\{\skal{\alpha,q}-m+1,\dots,\skal{\alpha,q}\}&\quad\text{if }r\leq\op{ht}(\beta),\\
\{\skal{\alpha,q}-m,\dots,\skal{\alpha,q}\}&\quad\text{if }0<\op{ht}(\beta)\leq r,\\
\{\skal{\alpha,q}-m,\dots,\skal{\alpha,q}-1\}&\quad\text{if }r-h\leq\op{ht}(\beta)<0,\\
\{\skal{\alpha,q}-m-1,\dots,\skal{\alpha,q}-1\}&\quad\text{if }\op{ht}(\beta)\leq r-h.
\end{cases}
\end{align*}
The claim now follows from \refl{kalpha}.
Note that in the last two cases $\skal{\alpha,q}-1\geq0$ is ensured by \refl{dominant}.
\end{proof}

With this description of $t_{\alpha}^p(\omega)$ we are able to verify that the rational Shi tableau defined in \refq{talpha} really specialises to the rational Shi tableau defined in \refq{tij} for type $A_{n-1}$.

\begin{myprop}{tijalphaA} Let $n,p$ be positive coprime integers, $i,j\in[n]$ with $i<j$, and $\omega\in\affS_n^p$ be a dominant $p$-stable affine permutation. Then $t_{i,j}^p(\omega)=t_{e_i-e_j}^p(\omega)$, which counts the number of affine inversions $(a,b)$ of $\omega$ such that $b<a+p$, and $\omega(a)\equiv i$ and $\omega(b)\equiv j$ modulo $n$.
\end{myprop}
\begin{proof}
This is a consequence of \refp{kalpha} and \refp{talpha}.
\end{proof}

The following is the main conjecture of this paper.

\begin{myconj}{shi} Let $\Phi$ be an irreducible crystallographic root system with affine Weyl group $\W$ and Coxeter number $h$, and let $p$ be a positive integer relatively prime to $h$. Then each dominant $p$-stable element of the affine Weyl group is determined uniquely by its rational Shi tableau.
\end{myconj}

\refcj{shi} is known in to be true in the Fu{\ss}--Catalan case where $p=mn+1$. In the next section we exploit the connections of the affine symmetric group to the combinatorics of rational Dyck paths to prove the conjecture if $\Phi$ is of type $A_{n-1}$. That is, we prove the following theorem. (The proof is found at the end of \refs{codinv}.)

\begin{mythrm}{shiA} Let $n,p$ be two positive coprime integers. Then each dominant $p$-stable affine permutation in $\affS_n$ is determined uniquely by its rational Shi tableau.
\end{mythrm}

We remark that the sum of the entries of the Shi tableau generalises the height statistic used by Stump in~\cite[Conj.~3.14]{Stump2010}.
Thus the $q$-Fu{\ss}--Catalan numbers proposed therein are generalised by the polynomials $q^{(p-1)r/2}C_{\Phi,p}(q^{-1})$, where $r$ is the rank of $\Phi$ and
\begin{align}
C_{\Phi,p}(q)
=\sum_{\omega} \prod_{\alpha\in\Phi^+}q^{t_{\alpha}^p(\omega)}
\end{align}
the sum being taken over all dominant $p$-stable elements of the affine Weyl group.

\smallskip
Pak and Stanley~\cite{Stanley1996,Stanley1998}
found a bijection between the regions of the ($m$-extended) Shi arrangement of type $A_{n-1}$ and the set of ($m$-)parking functions of length $n$. Gorsky, Mazin and Vazirani~\cite[Def.~3.8]{GMV2014} generalised this bijection to a map from $\affS_n^p$ to the set of rational parking functions.
Their \emph{rational Pak--Stanley labelling} $f(\omega)$ is defined by
\begin{align*}
f_i(\omega)
=\#\big\{(a,b)\in[n]\times\N:a<b<a+p,\omega(a)>\omega(b)\text{ and }\omega(b)\equiv i\text{ modulo }n\big\},
\end{align*}
where $i\in[n]$.
If $\omega$ is a dominant $p$-stable affine permutation then the Pak--Stanley labelling is obtained by taking the row-sums of the Shi tableau of $\omega$.
That is,
\begin{align*}
f_j(\omega)=\sum_{i=1}^{j-1}t_{i,j}^p(\omega).
\end{align*}
Consequently the \emph{dual Pak--Stanley labelling}, which we define by $f_i^*(\omega)=f_i(\omega^*)$, is obtained by taking the column-sums of the Shi tableau of $\omega$.
That is,
\begin{align*}
f_i^*(\omega)=\sum_{j=n-i+2}^nt_{n-i+1,j}^p(\omega).
\end{align*}
Partly the motivation for studying rational Shi tableaux comes from the fact that in view of the above identities they can be regarded as an intermediate step between dominant $p$-stable affine permutations and their images under the Pak--Stanley labelling.
If one replaces a dominant $p$-stable affine permutation $\omega$ by its rational Shi tableau then some information is apparently lost, since not all inversions of $\omega$ are taken into account.
Moving on to the Pak--Stanley labelling $f(\omega)$ gives up even more information on the nature of these inversions.
Nevertheless there seems to be just enough information left to determine $\omega$ uniquely.
While the injectivity of the Pak--Stanley labelling on the set of all $p$-stable affine permutations remains an open problem~\cite[Conj.~1.4]{GMV2014}, it follows from the work of Thomas and Williams~\cite{ThoWil2016} that the Pak--Stanley labelling is injective on the set of dominant $p$-stable affine permutations (see also the remarks at the end of \refs{codinv}).
In view of \reft{shiA} it is an interesting question whether Shi tableaux can offer new insights regarding this problem.

Another major open problem concerning the Pak--Stanley labelling is
to find an analogous labelling for different affine Weyl groups.
Notably also this problem can be solved for rational Shi tableaux (recall the definition in~\refq{talpha}).

\smallskip
To conclude this section we want to connect the world of cores and abaci to dominant affine permutations.
This connection is achieved by the next theorem, which follows essentially from the work of Lascoux~\cite{Lascoux2001}.
Let $\omega\in\affS_n$ be a dominant affine permutation, and define
$\gamma(\omega)=\{z\in\Z:\omega(z)\leq0\}$.

\begin{mythrm}{gamma} The map $\gamma$ is a bijection between dominant affine permutations in $\affS_n$ and balanced $n$-flush abaci.
\end{mythrm}

It is easy to adapt \reft{gamma} to involve $p$-stable affine permutations.

\begin{myprop}{pgamma} Let $\omega\in\affS_n$ be a dominant affine permutation. Then the abacus $\gamma(\omega)$ is $p$-flush if and only if $\omega\in\affS_n^p$ is $p$-stable.
\end{myprop}

Recall the map $\alpha$ from cores to balanced flush abaci that was introduced in~\refs{conj}. We obtain a bijection $\alpha^{-1}\circ\gamma$ between the dominant affine permutations in $\affS_n$ and $n$-cores. In particular $\alpha^{-1}\circ\gamma$ restricts to a bijection between dominant $p$-stable affine permutations and $\mathfrak C_{n,p}$. The following result appears to be new.

\begin{myprop}{conjcore} Let $\omega\in\affS_n$ be a dominant affine permutation. Then $\alpha^{-1}\circ\gamma(\omega^*)$ is the conjugate partition of $\alpha^{-1}\circ\gamma(\omega)$.
\end{myprop}

\begin{proof} This is best understood using a different description of the map $\alpha^{-1}\circ\gamma$. Following \cite[Sec.~1.2]{kschur} we read the one-line notation of $\omega$ from left to right, drawing a North step for each encountered non-positive number and an East step for each encountered positive number. The resulting path $P$ outlines the South-West boundary of the partition $\alpha^{-1}\circ\gamma(\omega)$.

By \refl{aut}~\refi{window}
$w^*(i)
=1-w(1-i)$
Hence $w^*(i)\leq0$ if and only if $w(1-i)$ is positive. Reading the one-line notation of $\omega^*$ from left to right and drawing a path as prescribed therefore yields the reverse path of $P$ with North and East step exchanged.
\end{proof}

For example, consider the dominant affine permutation $\omega=[-2,15,-1,16,-14,10,4]\in\affS_7^{16}$. 
Reading only the signs (zero counting as a minus) of the values of $\omega$ we obtain the sequence
\begin{align*}
\cdots-+-+----+-+-+-[-+-+-++]++++-++++++-+\cdots
\end{align*}
where the initial dots encode minuses, and the dots at the end encode pluses. Drawing the associated path we obtain the boundary of the $7,16$-core in \reff{core}.

\section{The codinv statistic}\label{Section:codinv}

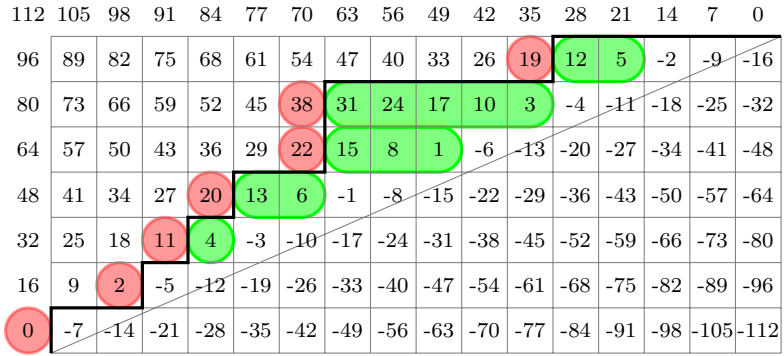
\begin{figure}[ht]
\begin{center}
\begin{tikzpicture}[scale=.6]
\draw[very thick, red, fill=red, opacity=.4, rounded corners=3mm]
(-1,0) rectangle (0,1)
(1,1) rectangle (2,2)
(2,2) rectangle (3,3)
(3,3) rectangle (4,4)
(5,4) rectangle (6,5)
(5,5) rectangle (6,6)
(10,6) rectangle (11,7)
;
\draw[very thick, green, fill=green!50, rounded corners=3mm]
(3,2) rectangle (4,3)
(4,3) rectangle (6,4)
(6,4) rectangle (9,5)
(6,5) rectangle (11,6)
(11,6) rectangle (13,7)
;
\draw[gray]
(0,0) grid (16,7)
(0,0)--(16,7)
;
\draw[very thick] (0,0)--(0,1)--(2,1)--(2,2)--(3,2)--(3,3)--(4,3)--(4,4)--(6,4)--(6,6)--(11,6)--(11,7)--(16,7)
;
\foreach \x in {0,1,2,3,4,5,6,7,8,9,10,11,12,13,14,15,16}{
\foreach \y in {0,1,2,3,4,5,6,7}{
\pgfmathtruncatemacro{\z}{16*\y - 7*\x}
\draw[xshift=-5mm,yshift=5mm] (\x,\y) node{\footnotesize{\z}};
}}
\end{tikzpicture}
\caption{A rational Dyck path $x\in\mathfrak D_{7,16}$.}
\label{Figure:dyck}
\end{center}
\end{figure}

Let $n$ and $p$ be positive coprime integers. A \emph{rational Dyck path} is a lattice path $x$ that starts at $(0,0)$, consists of $n$ North steps $N=(0,1)$ and $p$ East steps $E=(1,0)$, and never goes below the diagonal with rational slope $n/p$. Denote the set of all rational Dyck paths by $\mathfrak{D}_{n,p}$.

For $(i,j)$ with $0\leq i\leq p$ and $0\leq j\leq n$ place the label $\ell_{i,j}=jp-in$ in the unit square with bottom right corner $(i,j)$.
Given a rational Dyck path $x\in\mathfrak D_{n,p}$, we assign to each of its steps the label $\ell_{i,j}$, where $(i,j)$ is the starting point of the step.
Let $(\ell_1,\ell_2,\dots,\ell_n)$ be the vector consisting of the labels of the North steps of $x$ (indicated in red in \reff{dyck}) ordered increasingly.
Thus in \reff{dyck} we have $(\ell_1,\ell_2,\ell_3,\ell_4,\ell_5,\ell_6,\ell_7)=(0,2,11,19,20,22,38)$.
Denote by $H(x)$ the set of positive labels below $x$ (indicated in green in \reff{dyck}).
These correspond exactly to those boxes below $x$ but strictly above the diagonal. Hence we have $\op{area}(x)=\#H(x)$.

We introduce the following definitions. A \emph{codinv pair} of $x$ is a pair of integers $(\ell,h)$ such that $\ell$ is the label of a North step of $x$ and $h\in H(x)$ and $\ell\leq h\leq\ell+p$.
The \emph{codinv tableau} of $x$ is the collection of numbers
\begin{align}
d_{i,j}(x)
=\#\big\{(\ell_i,h):h\in H(x),\ell_i<h<\ell_i+p\text{ and }h\equiv\ell_j\mod n\big\}
\end{align}
where $i,j\in[n]$ with $i<j$.
For example consider the Dyck path in \reff{dyck}, which has codinv pairs
\begin{align*}
\begin{matrix}
(0,4),(0,6),(0,13),(0,1),(0,8),(0,15),(0,3),(0,10),(0,5),(0,12),\\
(2,4),(2,6),(2,13),(2,8),(2,15),(2,3),(2,10),(2,17),(2,5),(2,12),\\
(11,13),(11,15),(11,17),(11,24),(11,12),\\
(19,24),(19,31),(20,24),(20,31),(22,24),(22,31).
\end{matrix}
\end{align*}
Its codinv tableau is found in \reff{codinv}.

\begin{figure}[t]
\begin{center}
\begin{tikzpicture}[scale=.6]
\begin{scope}
\draw[xshift=5mm,yshift=5mm]
(0,0) node{\footnotesize{$d_{1,2}$}}
(0,1) node{\footnotesize{$d_{1,3}$}}
(0,2) node{\footnotesize{$d_{1,4}$}}
(0,3) node{\footnotesize{$d_{1,5}$}}
(0,4) node{\footnotesize{$d_{1,6}$}}
(0,5) node{\footnotesize{$d_{1,7}$}}
(1,1) node{\footnotesize{$d_{2,3}$}}
(1,2) node{\footnotesize{$d_{2,4}$}}
(1,3) node{\footnotesize{$d_{2,5}$}}
(1,4) node{\footnotesize{$d_{2,6}$}}
(1,5) node{\footnotesize{$d_{2,7}$}}
(2,2) node{\footnotesize{$d_{3,4}$}}
(2,3) node{\footnotesize{$d_{3,5}$}}
(2,4) node{\footnotesize{$d_{3,6}$}}
(2,5) node{\footnotesize{$d_{3,7}$}}
(3,3) node{\footnotesize{$d_{4,5}$}}
(3,4) node{\footnotesize{$d_{4,6}$}}
(3,5) node{\footnotesize{$d_{4,7}$}}
(4,4) node{\footnotesize{$d_{5,6}$}}
(4,5) node{\footnotesize{$d_{5,7}$}}
(5,5) node{\footnotesize{$d_{6,7}$}}
;
\end{scope}
\draw(7,3)node{$=$};
\begin{scope}[xshift=8cm]
\draw[xshift=5mm,yshift=5mm]
(0,0) node{\footnotesize{$0$}}
(0,1) node{\footnotesize{$1$}}
(0,2) node{\footnotesize{$2$}}
(0,3) node{\footnotesize{$2$}}
(0,4) node{\footnotesize{$3$}}
(0,5) node{\footnotesize{$2$}}
(1,1) node{\footnotesize{$1$}}
(1,2) node{\footnotesize{$2$}}
(1,3) node{\footnotesize{$2$}}
(1,4) node{\footnotesize{$2$}}
(1,5) node{\footnotesize{$3$}}
(2,2) node{\footnotesize{$1$}}
(2,3) node{\footnotesize{$1$}}
(2,4) node{\footnotesize{$1$}}
(2,5) node{\footnotesize{$2$}}
(3,3) node{\footnotesize{$0$}}
(3,4) node{\footnotesize{$0$}}
(3,5) node{\footnotesize{$2$}}
(4,4) node{\footnotesize{$0$}}
(4,5) node{\footnotesize{$2$}}
(5,5) node{\footnotesize{$2$}}
;
\end{scope}
\end{tikzpicture}
\caption{The codinv tableau of the Dyck path in \reff{dyck}.}
\label{Figure:codinv}
\end{center}
\end{figure}

\sk
Similar constructions have appeared before.
First we remark that the codinv tableau is related to, albeit not the same as, the laser fillings of Ceballos, Denton and Hanusa~\cite[Def.~5.13]{CDH2015}. The row-sums and column-sums of the codinv tableau and the laser filling of a Dyck path agree.
However, the codinv tableau is always of staircase shape while the laser filling sits inside the boxes below the rational Dyck path.
Secondly we note that codinv pairs have been considered by Gorsky and Mazin~\cite{GM2013} using slightly different notation.
However, they only considered the column-sums of the codinv tableau.
See also the remarks following \reft{zeta}.

\smallskip
Our next aim is to relate rational Shi tableaux to codinv tableaux.
Our starting point is an  elegant bijection $\varphi:\mathfrak D_{n,p}\to\mathfrak C_{n,p}$
between rational Dyck paths and simultaneous cores discovered by Anderson~\cite{Anderson2002}.
For any finite set $H=\{h_1,\dots,h_k\}$ of positive integers there exists a unique partition $\lambda$ such that $H$ is the set of hook lengths of the cells in the first column of $\lambda$. For any $x\in\mathfrak D_{n,p}$ let $\varphi(x)$ be the partition such that the set of hook lengths of the cells in its first column equals $H(x)$. The Dyck path in \reff{dyck} thereby corresponds to the core in \reff{core}.

\begin{mythrm}{anderson}~\textnormal{\cite[Prop.~1]{Anderson2002}} The map $\varphi:\mathfrak D_{n,p}\to\mathfrak C_{n,p}$ is a bijection.
\end{mythrm}

Gorsky, Mazin and Vazirani~\cite[Sec.~3.1]{GMV2014} defined a generalised Anderson map $\mathcal A$ 
that bijectively maps the set of $p$-stable affine permutations to the set of rational parking functions. We are interested in the restriction of the Anderson map to the set of dominant $p$-stable affine permutations, which can be written as the composition $\mathcal A=\varphi^{-1}\circ\alpha^{-1}\circ\gamma$ of maps we have already discussed. Here the inverse of Anderson's bijection $\varphi$ makes an appearance, which accounts for the name of the Anderson map.
We use the Anderson map to relate rational Shi tableaux to the codinv statistic.

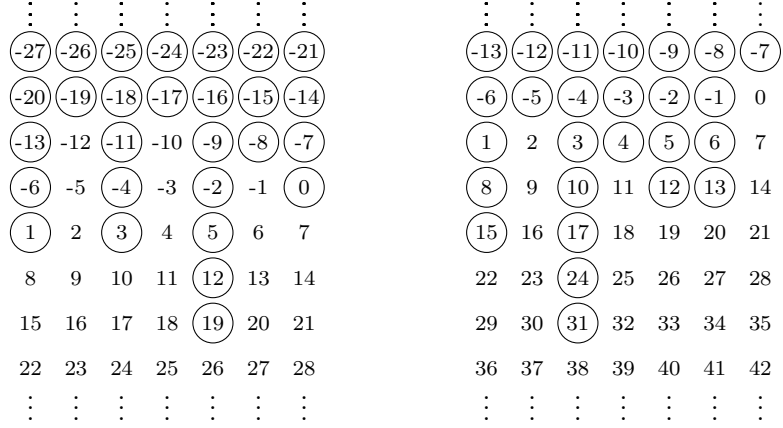
\begin{figure}[t]
\begin{center}
\begin{tikzpicture}[scale=.6]
\pgfmathtruncatemacro{\n}{7}
\begin{scope}[]
\pgfmathtruncatemacro{\min}{4}
\pgfmathtruncatemacro{\max}{-3}
\drawrunner{1}{0}
\drawrunner{2}{-3}
\drawrunner{3}{0}
\drawrunner{4}{-3}
\drawrunner{5}{2}
\drawrunner{6}{-2}
\drawrunner{7}{-1}
\end{scope}
\begin{scope}[xshift=10cm,yshift=2cm]
\pgfmathsetmacro{\min}{2}
\pgfmathsetmacro{\max}{-5}
\drawrunner{1}{2}
\drawrunner{2}{-1}
\drawrunner{3}{4}
\drawrunner{4}{0}
\drawrunner{5}{1}
\drawrunner{6}{1}
\drawrunner{7}{-2}
\end{scope}
\end{tikzpicture}
\caption{The balanced abacus $A=\gamma(\omega)$, where $\omega=[-2,15,-1,16,-14,10,4]$, (left) and the normalised abacus $B=\beta\circ\alpha^{-1}(A)$ (right), both depicted on $7$ runners.
With the notation of the proof of \reft{dinvshi} we have $\omega^{-1}=[-12,-10,-1,7,8,10,26]$ and $s=[2,4,6,7,1,3,5]$ and $\sigma=[7,2,4,5,6,1,3]$.
Moreover $\varphi^{-1}\circ\alpha^{-1}(A)=\varphi^{-1}\circ\beta^{-1}(B)=\mathcal{A}(\omega)$ is the rational Dyck path of \reff{dyck}.}
\label{Figure:abaci}
\end{center}
\end{figure}

\begin{mythrm}{dinvshi} Let $n,p$ be positive coprime integers and $\omega\in\affS_n^p$ be a dominant $p$-stable affine permutation. Then the rational Shi tableau of $\omega$ equals the codinv tableau of $\mathcal A(\omega)$. That is,
$t_{i,j}^p(\omega)=d_{i,j}(\mathcal A(\omega))$
for all $i,j\in[n]$ with $i<j$.
\end{mythrm}

\begin{proof} Let $\omega^{-1}=t_qs$ where $q\in\coQ$ and $s\in\S_n$, and fix $i,j\in[n]$ such that $i<j$. By \refp{tijalphaA} $t_{i,j}^p(\omega)$ equals the number of inversions $(j,kn+i)$ of $\omega^{-1}$ such that $\omega^{-1}(j)-\omega^{-1}(kn+i)<p$.

Let $A=\gamma(\omega)$ be an abacus on $n$ runners (see \reff{abaci}).
Then $\omega^{-1}(j)$ is the minimal gap of $A$ in the runner $s(j)$.
Moreover, $(j,kn+i)$ is an inversion of $\omega^{-1}$ contributing to $t_{i,j}^p(\omega)$ if and only if $\omega^{-1}(kn+i)$ is a non-minimal gap of $A$ in the runner $s(i)$ and
\begin{align*}
\omega^{-1}(j)-p<\omega^{-1}(nk+i)<\omega^{-1}(j).
\end{align*}
Hence $t_{i,j}^p(\omega)$ counts the number of non-minimal gaps $g$ in runner $s(i)$ such that $m-p<g<m$ where $m$ is the minimal gap in runner $s(j)$.

Equivalently $t_{i,j}^p(\omega)$ counts the number of beads $b$ in runner $s(j)$ such that $m<b<m+p$ where $m$ is the minimal gap in runner $s(i)$.
Define another abacus on $n$ runners, namely
\begin{align*}
B=\beta\circ\alpha^{-1}\circ\gamma=\{z+\ell-1:z\in A\},
\end{align*}
where $\ell$ is the length of the partition $\alpha^{-1}\circ\gamma(\omega)$. Moreover define $\sigma\in\S_n$ by $\sigma(i)\equiv s(i)+\ell-1$ modulo $n$. Then $t_{i,j}^p(\omega)$ counts the number of beads $b$ in the runner $\sigma(j)$ of $B$ such that $m<b<m+p$ where $m$ is the minimal gap in the runner $\sigma(i)$ of $B$.

Since $B$ is normalised, the minimal gap of each runner of $B$ is non-negative.
Thus it is the same to consider only positive beads of $B$.
But the positive beads of $B$ are just the hook-lengths of the cells in the first column of $\alpha^{-1}\circ\gamma$ and therefore make up the set $H(\mathcal A(\omega))$. Moreover the minimal gaps of the runners of $B$ are just the labels of the North steps of $\mathcal A(\omega)$.

The theorem now follows from the observation that $\sigma$ sorts the minimal gaps of $B$ increasingly. This is implied by the fact that $s$ sorts the minimal gaps of $A$ increasingly, which form the window of the affine Gra{\ss}mannian permutation $\omega^{-1}$.
\end{proof}

From \reft{dinvshi}, \refp{shiconj} and \refp{conjcore} we obtain an interesting result for free.

\begin{mycor}{conjdinv} Let $n,p$ be positive coprime integers and $x\in\mathfrak{D}_{n,p}$ be a rational Dyck path. Then the codinv tableau of $x$ is the transpose of the codinv tableau of $\varphi^{-1}(\varphi(x)')$.
\end{mycor}

Note that the map $x\mapsto\varphi^{-1}(\varphi(x)')$, which is evidently an involution on $\mathfrak{D}_{n,p}$, has been called the \emph{rank compliment} by
Xin 
while Ceballos, Denton and Hanusa 
call it \emph{conjugation} on Dyck paths.


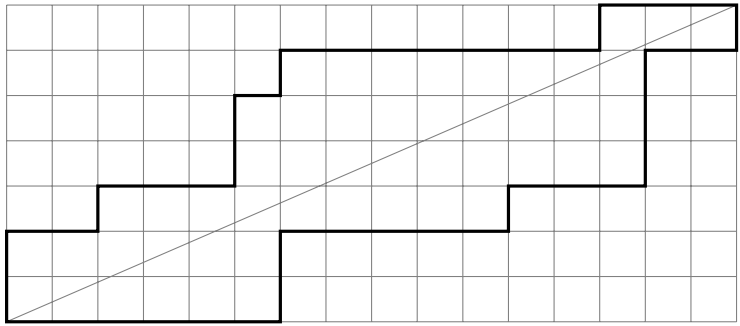
\begin{figure}[t]
\begin{center}
\begin{tikzpicture}[scale=.6]
;
\draw[gray]
(0,0) grid (16,7)
(0,0)--(16,7)
;
\draw[very thick] (0,0)--(0,2)--(2,2)--(2,3)--(5,3)--(5,5)--(6,5)--(6,6)--(13,6)--(13,7)--(16,7)-
-(16,6)--(14,6)--(14,3)--(11,3)--(11,2)--(6,2)--(6,0)--cycle
;
%
\end{tikzpicture}
\caption{The image $\zeta(x)$ of the rational Dyck path in \reff{dyck} under the zeta map, and the (rotated) path $\eta(x)$ below the diagonal.}
\label{Figure:zeta}
\end{center}
\end{figure}

\smallskip
The zeta map $\zeta:\mathfrak D_{n,p}\to\mathfrak D_{n,p}$ is a map on rational Dyck paths that has appeared many times in the literature. Andrews, Krattenthaler, Orsina and Papi~\cite{AKOP2002} first defined the inverse of the zeta function in the Catalan case $\mathfrak D_{n,n+1}$. Haglund~\cite{Haglund2003} 
studied the Catalan instance of the zeta map in connection with diagonal harmonics. More recently multiple descriptions of the rational zeta map have been found (see~\cite{GM2013,ALW2014:sweep_maps}).

The \emph{zeta map} can be defined as follows.
Let $x=x_1x_2\dots x_{n+p}\in\mathfrak{D}_{n,p}$ be a rational Dyck path with steps $x_i\in\{N,E\}$. For $i\in[n+p]$ let $\ell(x_i)$ be the label of the step $x_i$.
Then
$\zeta(x)=x_{\sigma(1)}x_{\sigma(2)}\dots x_{\sigma(n+p)}$,
where $\sigma\in\S_{n+p}$ is the unique permutation such that $\ell(x_{\sigma(1)})<\ell(x_{\sigma(2)})<\dots<\ell(x_{\sigma(n+p)})$.
The zeta map is accompanied by a second map on rational Dyck paths that is called the \emph{eta map} in~\cite{CDH2015}. It is defined by
$\eta(x)=\zeta\circ\varphi^{-1}((\varphi(x)')$
where $\lambda'$ is the conjugate partition of $\lambda$ and $\varphi$ is Anderson's bijection.
See \reff{zeta}.

\smallskip
Using codinv tableaux we reprove a connection between the Anderson map, the zeta map and the Pak--Stanley labelling that was already observed in \cite[Thm.~5.3]{GMV2014}.
Note that each rational Dyck path $x\in\mathfrak D_{n,p}$ is the South-East boundary of a partition $\lambda$ fitting inside the $p\times n$ rectangle. We call $\lambda$ the \emph{complement} of $x$. For example, the complement of the Dyck path $x$ in \reff{dyck} is the partition $\lambda=(11,6,6,4,3,2,0)$.
The complements of the Dyck paths $\zeta(x)$ and $\eta(x)$ in \reff{zeta} are $(13,6,5,5,2,0,0)$ and $(10,10,5,2,2,2,0)$ respectively.

\begin{mythrm}{zeta} Let $n,p$ be positive coprime integers and $\omega\in\affS_n^p$ be a dominant $p$-stable affine permutation. Then the complement of $\zeta\circ\mathcal A(\omega)$ equals the (reversed\footnote{In our convention partitions are weakly decreasing while $f(\omega)$ is weakly increasing.}) Pak--Stanley labelling $f(\omega)$. Moreover, the complement of $\eta\circ\mathcal A(\omega)$ equals the (reversed) dual Pak--Stanley labelling $f^*(\omega)$.
\end{mythrm}

\begin{proof} Let $x=\mathcal{A}(\omega)$.
By \reft{dinvshi} the Pak--Stanley labelling $f(\omega)$ is given by the row-sums of the codinv tableau of $x$.

Consider the $j$-th North step of $\zeta(x)$. An East step of $x$ precedes this North step in $\zeta(x)$ if and only if its label is less than $\ell_j$.
Let $L$ be the set of lines of slope $n/p$ going through the bottom right corner of a cell labelled by $\ell_j-kn$ for some $k>0$.
The label of an East step of $x$ is less than $\ell_j$ if and only if the East step is intersected by a line in $L$.

Moreover each East step of $x$ is intersected by at most one such line.
Thus the number of East steps preceding the $j$-th North step of $\zeta(x)$ is counted by the number of intersections of an East step of $x$ and a line in $L$.

For each line in $L$ the number of intersected East steps equals the number of intersected North steps.
Thus the number of East steps preceding the $j$-th North step of $\zeta(x)$ is counted by the number of intersections of a North step of $x$ and a line in $L$.

It takes a moment of thought to verify that this number is given by the sum $\sum_{i=1}^{j-1}d_{i,j}(x)$ of a row of the codinv tableau.
Indeed the number $d_{i,j}(x)$ counts the number of intersections of the North step of $x$ labelled $\ell_i$ and a line in $L$.
We conclude that the $j$-th part of the complement of $\zeta(x)$
equals $f_{n-j+1}(\omega)$.

\smallskip
The dual statement concerning the eta map could be deduced analogously. However, it follows immediately using the involutive automorphism (see \refp{shiconj} and \refp{conjcore}).
\end{proof}

Our proof is inspired by the resemblance of codinv tableaux and laser fillings.
\reft{zeta} should therefore be compared to~\cite[Thm.~5.15]{CDH2015}.
Moreover we remark that taking the column-sums of the codinv tableau actually coincides with the definition of a map of Gorsky and Mazin~\cite[Def.~3.3]{GM2013}, which is the eta map in our notation.
%
Thus it is well worth comparing \reft{zeta} to~\cite[Thm.~4.12]{ALW2014:sweep_maps}.

\smallskip
It is well known that $\op{skl}(\varphi(x))=(n-1)(p-1)/2-\op{area}(\zeta(x))$. This follows from Armstrong's definition of the zeta map~\cite[Sec.~4.2]{ALW2014:sweep_maps}, which uses cores.
As a first consequence of \reft{zeta} we obtain a justification for calling the codinv tableau a refinement of the skew length statistic.

\begin{mycor}{codinv} Let $x\in\mathfrak{D}_{n,p}$. Then
\begin{align*}
\sum_{i,j}d_{i,j}(x)&=\frac{(n-1)(p-1)}{2}-\op{area}(\zeta(x)).
\end{align*}
\end{mycor}

Furthermore we may now use a result of Ceballos, Denton and Hanusa 
to prove that each rational Dyck path is determined uniquely by its codinv tableau, and equivalently, that each dominant $p$-stable affine permutation is determined uniquely by its rational Shi tableau.

\begin{mythrm}{unique} 
Let $x,y\in\mathfrak D_{n,p}$ 
such that $d_{i,j}(x)=d_{i,j}(y)$ for all $i,j\in[n]$ with $i<j$. Then $x=y$.
\end{mythrm}

\begin{proof}[Proof of Theorems \ref{Theorem:shiA} and \ref{Theorem:unique}] We deduce the claims from \cite[Thm.~6.3]{CDH2015}, which asserts that any rational Dyck path $x$ can be reconstructed from the pair $(\zeta(x),\eta(x))$. By Theorems~\ref{Theorem:dinvshi} and~\ref{Theorem:zeta} the codinv tableau of $x$ encodes both $\zeta(x)$ and $\eta(x)$ in terms of column-sums and row-sums. Therefore it contains enough information to determine the path $x$ uniquely. Using the Anderson map again we obtain the analogous statement for rational Shi tableaux and dominant $p$-stable affine permutations.
\end{proof}

We close with an exciting recent result due to Thomas and Williams.

\begin{mythrm}{zetabij}~\textnormal{\cite[Cor.~6.4]{ThoWil2016}} The zeta map $\zeta:\mathfrak{D}_{n,p}\to\mathfrak{D}_{n,p}$ is a bijection.
\end{mythrm}

It should be clear from the above discussion that \reft{zetabij} implies that each dominant $p$-stable affine permutation is determined uniquely by its Pak--Stanley labelling. Thus the stronger \reft{zetabij} could replace \cite[Thm.~6.3]{CDH2015} in the proof of \reft{shiA}.

\section{Acknowledgements}

The author is grateful to Cesar Ceballos, Jim Haglund, Christian Krattenthaler and Marko Thiel for many insightful discussions.



\bibliographystyle{alpha}
\bibliography{/users/rsulzg/Documents/Latex/robinbib}

\newcommand{\etalchar}[1]{$^{#1}$}
\begin{thebibliography}{{Hum}90}

\bibitem[AHJ14]{AHJ2014}
Drew {Armstrong}, Christopher~R.~H. {Hanusa}, and Brant~C. {Jones}.
\newblock Results and conjectures on simultaneous core partitions.
\newblock {\em European J. Combin.~\myfontB{41}}, pages 205--220, 2014.
\newblock \texttt{arXiv:1308.0572}.

\bibitem[AKOP02]{AKOP2002}
George {Andrews}, Christian {Krattenthaler}, Luigi {Orsina}, and Paolo {Papi}.
\newblock {$ad$}-nilpotent {$\mathfrak{b}$}-ideals in {$sl(n)$} having fixed
  class of nilpotence: {Combinatorics} and enumeration.
\newblock {\em Trans. Amer. Math. Soc.~\myfontB{354}}, pages 3835--3853, 2002.

\bibitem[ALW15]{ALW2014:sweep_maps}
Drew {Armstrong}, Nicholas~A. {Loehr}, and Gregory~S. {Warrington}.
\newblock Sweep maps: A continuous family of sorting algorithms.
\newblock {\em Adv. Math.~\myfontB{284}}, pages 159--185, 2015.

\bibitem[{And}02]{Anderson2002}
Jaclyn {Anderson}.
\newblock Partitions which are simultaneously {$t_1$}- and {$t_2$}-cores.
\newblock {\em Discrete Math.~\myfontB{248}}, pages 237--243, 2002.

\bibitem[BB96]{BjoBre1996}
Anders {Bj\"orner} and Francesco {Brenti}.
\newblock {Affine permutations of type~{$A$}}.
\newblock {\em Electron. J. Combin.~\myfontB[\,(2)]{3}, The Foata Festschrift
  volume}, page R18, 1996.

\bibitem[BB05]{BjoBre}
Anders {Bj\"orner} and Francesco {Brenti}.
\newblock {\em Combinatorics of {Coxeter} Groups}, volume 231 of {\em Grad.
  Texts in Math.}
\newblock Springer, New York, 2005.

\bibitem[CDH16]{CDH2015}
Cesar {Ceballos}, Tom {Denton}, and Christopher~R.~H. {Hanusa}.
\newblock Combinatorics of the zeta map on rational {Dyck} paths.
\newblock {\em J.~Combin. Theory Ser.~A}, 141:33--77, 2016.
\newblock \texttt{arXiv:1504.06383v3}.

\bibitem[FTV11]{FTV2011}
Susanna {Fishel}, Eleni {Tzanaki}, and Monica {Vazirani}.
\newblock Counting {Shi} regions with a fixed separating wall.
\newblock In {\em Proceedings of the 23rd International Conference on Formal
  Power Series and Algebraic Combinatorics (FPSAC), Reykjav{\'i}k, Iceland},
  pages 351--362, 2011.

\bibitem[GM13]{GM2013}
Eugene {Gorsky} and Mikhail {Mazin}.
\newblock Compactified {Jacobians} and {$q,t$}-{Catalan} numbers, {I}.
\newblock {\em J.~Combin. Theory Ser.~A~\myfontB{120}}, pages 49--63, 2013.

\bibitem[GMV16]{GMV2014}
Eugene {Gorsky}, Mikhail {Mazin}, and Monica {Vazirani}.
\newblock Affine permutations and rational slope parking functions.
\newblock {\em Trans. Amer. Math. Soc.}, 2016.
\newblock \texttt{arXiv:1403.0303}.

\bibitem[{Hag}03]{Haglund2003}
James {Haglund}.
\newblock Conjectured statistics for the $q,t$-{Catalan} numbers.
\newblock {\em Adv. Math.~\myfontB{175}}, pages 319--334, 2003.

\bibitem[{Hum}90]{Humphreys}
James~E. {Humphreys}.
\newblock {\em Reflection Groups and {Coxeter} Groups}, volume~29 of {\em
  Cambridge Stud. Adv. Math.}
\newblock Cambridge Univ. Press, Cambridge, 1990.

\bibitem[JK81]{JamesKerber}
Gordon {James} and Adalbert {Kerber}.
\newblock {\em The Representation Theory of the Symmetric Group}, volume~16 of
  {\em Encyclopedia Math. Appl.}
\newblock Addison--Wesley, Reading MA, 1981.

\bibitem[{Las}01]{Lascoux2001}
Alain {Lascoux}.
\newblock Ordering the affine symmetric group.
\newblock In {\em Proceedings of the Euroconference, Algebraic Combinatorics
  and Applications (ALCOMA), G{\"o}{\ss}weinstein, Germany, 1999}, pages
  219--231. Springer, Berlin, 2001.

\bibitem[LLM{\etalchar{+}}14]{kschur}
Thomas {Lam}, Luc {Lapointe}, Jennifer {Morse}, Anne {Schilling}, Mark
  {Shimozono}, and Mike {Zabrocki}.
\newblock {\em {$k$}-{Schur} Functions and Affine {Schubert} Calculus}.
\newblock Fields Inst. Monogr.~\myfontB{33}. Springer, New York, 2014.

\bibitem[{Shi}86]{Shi:Kazhdan_Lusztig}
Jian{-}Yi {Shi}.
\newblock {\em The {Kazhdan-Lusztig} Cells in Certain Affine {Weyl} Groups},
  volume 1179 of {\em Lecture Notes in Math.}
\newblock Springer, Berlin Heidelberg, 1986.

\bibitem[{Shi}87]{Shi1987:alcoves}
Jian{-}Yi {Shi}.
\newblock Alcoves types corresponding to an affine {Weyl} group.
\newblock {\em J. Lond. Math. Soc.~(2)~\myfontB[\,(1)]{35}}, pages 42--55,
  1987.

\bibitem[{Sta}96]{Stanley1996}
Richard~P. {Stanley}.
\newblock Hyperplane arrangements, interval orders, and trees.
\newblock {\em Proc. Nat. Acad. Sci.~\myfontB{93}}, pages 2620--2625, 1996.

\bibitem[{Sta}98]{Stanley1998}
Richard~P. {Stanley}.
\newblock Hyperplane arrangements, parking functions and tree inversions.
\newblock {\em Mathematical Essays in Honor of Gian-Carlo Rota}, pages
  359--375, 1998.

\bibitem[{Stu}10]{Stump2010}
Christian {Stump}.
\newblock {$q,t$}-{Fu\ss}--{Catalan} numbers for finite reflection groups.
\newblock {\em J.~Algebraic Combin.~\myfontB[\,(1)]{32}}, pages 67--97, 2010.

\bibitem[{Sul}16]{Sul2016}
Robin {Sulzgruber}.
\newblock Rational {Shi} tableaux and the skew length statistic.
\newblock In {\em DMTCS Proceedings of the 28th International Conference on
  Formal Power Series and Algebraic Combinatorics (FPSAC), Vancouver, Canada},
  pages 1147--1158, 2016.

\bibitem[{Thi}16]{Thiel2015}
Marko {Thiel}.
\newblock From {Anderson} to {Zeta}.
\newblock {\em Adv. in Appl. Math.~\myfontB{81}}, pages 156--201, 2016.
\newblock \texttt{arXiv:1504.07363v3}.

\bibitem[TW16]{ThoWil2016}
Hugh {Thomas} and Nathan {Williams}.
\newblock Sweeping up zeta.
\newblock 2016.
\newblock Preprint at \texttt{arXiv:1512.01483v4}.

\bibitem[{Xin}15]{Xin2015}
Guoce {Xin}.
\newblock Rank complement of rational {Dyck} paths and conjugation of
  {$(m,n)$}-core partitions.
\newblock 2015.
\newblock Preprint at \texttt{arXiv:1504.02075v2}.

\end{thebibliography}



\end{document}